%% file: main.tex
\documentclass[12pt,dvipsnames, DIV=13]{scrartcl}
\usepackage{amsmath,amsthm}
\usepackage[colorlinks,
urlcolor=Blue,
citecolor=JungleGreen,
linkcolor=Magenta]{hyperref}
\usepackage[algo2e,algosection,vlined,linesnumbered,resetcount,algoruled]{algorithm2e}
\usepackage[capitalise,noabbrev]{cleveref}
\newtheorem{theorem}{Theorem}[section]

\theoremstyle{definition}

\theoremstyle{remark}

\usepackage{amssymb,mathtools,subdepth}
\usepackage{tabularx,booktabs,multirow}
\newcolumntype{Y}{>{\centering\arraybackslash}X}

\usepackage{listings}
\usepackage[hang,flushmargin]{footmisc}
\makeatletter
\newcommand{\algorithmfootnote}[2][\footnotesize]{
  \let\old@algocf@finish\@algocf@finish
  \def\@algocf@finish{\old@algocf@finish
    \leavevmode\rlap{\begin{minipage}{\linewidth}
    #1#2
    \end{minipage}}
  }
}
\makeatother

\usepackage{subfigure}
\makeatletter
\renewcommand\p@subfigure{\thefigure}
\makeatother

\usepackage{pgfplots}
\pgfplotsset{compat=1.16}
\input{tikz.tex}
\input{tikz_lists.tex}

\usepackage[active]{srcltx}
\usepackage{MnSymbol}

\usepackage{tikz}
\usetikzlibrary{patterns}

\input{metadata.tex}

\input{commands.tex}

\newcommand{\funding}[1]{#1}
\newcommand{\email}[1]{\href{mailto:#1}{#1}}
\renewenvironment{abstract}{{\small \noindent\textbf{Abstract.}}}{\smallskip}
\newenvironment{keywords}{{\small \noindent\textbf{Keywords.}}}{\smallskip}
\newenvironment{MSCcodes}{{\small \noindent\textbf{MSC codes.}}}{\smallskip}
\date{\vskip -30pt}

\title{\thetitle\thanks{\thetitlenote{} \funding{\thefunding}}}
\author{\theauthori\thanks{\theaffiliationi{} (\email{\theemaili}).} \and
\theauthorii\thanks{\theaffiliationii{} (\email{\theemailii}).} \and
\theauthoriii\thanks{\theaffiliationiii{}.
}
\and
\theauthoriv\thanks{\theaffiliationiv{} (\email{\theemailiv}).}}

\begin{document}
  \maketitle

  \begin{abstract}
    \theabstract
  \end{abstract}



  \section{Introduction}
  A Sylvester matrix equation has the form
  \begin{equation}
    \label{eq:sylv}
    AX + XB = C,
  \end{equation}
  where the coefficients $A \in \Cmm$ and $B \in \Cnn$ and the
  right-hand side $C \in \Cmn$ are given, whilst $X \in \Cmn$ is the unknown.
  The matrix equation~\eqref{eq:sylv} can be recast as the $mn \times
  mn$ linear system
  \begin{equation}
    \label{eq:sylv-ls}
    M_f \vvec(X) = \vvec(C),\quad
    M_f := I_n \otimes A + B^T \otimes I_m,
  \end{equation}
  where $I_k \in \Ckk$ denotes the identity matrix of order $k$,
  $\otimes$ is the infix operator that computes the Kronecker product,
  and $\vvec$ is the operator that stacks the columns of an $m \times n$
  matrix into a vector of length~$mn$.
  The subscript $f$ stands for ``full'' and stresses that the Kronecker matrix in~\eqref{eq:sylv-ls} corresponds to the full Sylvester equation---in later sections we will use the subscript $t$ for ``triangular'' to denote the Kronecker matrix of a triangular equation.
A special case of~\eqref{eq:sylv} is the continuous Lyapunov equation, which has the form
  \begin{equation}
  	\label{eq:lyap}
  	AX + XA^{*} = C.
  \end{equation}

  The matrix equations~\eqref{eq:sylv} and~\eqref{eq:lyap} have both been extensively investigated in the literature, and theoretical results, as well as algorithms and software tools, are available.
  This is motivated by the key role that these two equations play in both theory and applications.
  Theoretically, they are essential in block diagonalization~\cite[sect.~7.6.3]{gova13} and
  perturbation theory for matrix functions~\cite[sects.~5.1 and 6.1]{high:FM}.
  In applications, they are employed in signal
  processing~\cite{care96},~\cite{vand91},
  system balancing~\cite{lhpw87},
  control~\cite{cofr03},~\cite{datt94},
  model reduction~\cite{babe08},~\cite{lhpw87},
  \cite{svr08},~\cite{soan02},
  machine learning~\cite{cswz08},
  numerical methods for matrix functions~\cite{dahi03},~\cite{hili21},
  and matrix differential Riccati equations~\cite{beme13}.
  We refer the reader to the survey by Bhatia and Rosenthal~\cite{bhro97} for a historical perspective and a review of theoretical results, and to the work of Simoncini~\cite{simo16} for a discussion on computational methods.

  Given the availability of new mixed-precision hardware, there has been a
  renewed interest in mixed-precision algorithms.
  Such algorithms have been studied and developed for a number of applications.
  For a discussion of mixed-precision methods for numerical linear algebra problems, we refer the reader to the surveys by Abdelfattah et al.~\cite{aabc21} and by Higham and Mary~\cite{hima22}.
  In this work, we design mixed-precision numerical methods for the solution
  of~\eqref{eq:sylv} and~\eqref{eq:lyap} that are as accurate as, but
  faster than, existing alternatives that rely only on one level of precision.
  The new algorithm are based on the Schur factorization and, therefore, they are mainly of interest when solving equations with dense coefficients of moderate size.
  We focus on the mixed-precision framework where
  \textit{two} floating-point arithmetics are involved.

  In the next section, we briefly recall the relevant preliminary results on these matrix equations and summarize the classical algorithms used to solve them.
  In \cref{sec:pert-sylv-algs}, we revisit a stationary iteration for linear systems, design an iterative refinement scheme for perturbed quasi-triangular Sylvester equations, and investigate its convergence conditions and attainable residual.
  In~\cref{sec:mixed-precision}, we use this scheme as a building block to construct mixed-precision algorithms for the Sylvester equation, and we study their convergence, computational complexity, and storage requirements.
  In~\cref{sec:flop-based-comp-model}, we establish a flop-based computational model to compare the cost of the new mixed-precision algorithms with that of the standard Bartels--Stewart algorithm in high precision.
  In~\cref{sec:gmres-ir}, we discuss the possibility of using GMRES-IR for solving the Sylvester equation in mixed precision.
  In~\cref{sec:numeric.experiments}, we test the numerical stability and the performance of the new mixed-precision algorithms on various Sylvester and Lyapunov equations from the literature.
  Conclusions are drawn in~\cref{sec:conclusions}.

  \section{Background}
  \label{sec:background}
  In this section, we introduce our notation and recall some background material needed later on.

  We use the standard model of floating-point arithmetic~\cite[sect.~2.2]{high:ASNA2}, and we consider two floating-point arithmetics with unit roundoffs $\ulow$ and $\uhigh$ that satisfy
  \begin{equation}
      \label{eq:prec-inequality}
      0 < \uhigh < \ulow < 1.
  \end{equation}
  We call these two arithmetics low and high precisions, respectively, and we say that we are ``computing in precision $\uhigh$'' to mean that $\uhigh$ is the unit roundoff of the current (working) precision.
  The features of the formats we consider are shown in \cref{tab:fp-parameters},
  where we list the number of binary digits in the significand (including the implicit leading bit), $t$, and in the exponent.
  The former directly affects the unit roundoff, $u=2^{-t}$,
  the maximum relative error introduced by rounding to nearest representable floating-point number, while the latter governs the dynamic range, which grows doubly exponentially in the number of exponent bits.

  \begin{table}[t]
    \caption{Parameters of five floating-point formats. The values of unit roundoff are to two significant digits, and approximate ranges of representable floating-point numbers are shown in the last column.}
    \label{tab:fp-parameters}
    \centering
    \begin{tabularx}{\linewidth}{@{}Xcccc@{}}
      \toprule
      \multirow{2}[3]{*}{Format} & \multicolumn{2}{c}{Number of bits}
      &  & \\
      \cmidrule{2-3}
    &  in significand $(t)$ &  in exponent  &  Unit roundoff $(u=2^{-t})$ &  Range \\
    \midrule
    bfloat16  &  \phantom{1}8  &   \phantom{1}8   &  $3.9\times 10^{-3\phantom{0}}$ &  $10^{\pm   38\phantom{0}}$ \\
    binary16  &   11  &   \phantom{1}5  & $4.9\times 10^{-4\phantom{0}}$  & $10^{\pm    5\phantom{00}}$  \\
    TensorFloat-32 &   11  &  \phantom{1}8   &  $4.9\times 10^{-4\phantom{0}}$ &  $10^{\pm    38\phantom{0}}$ \\
    binary32  &   24  &   \phantom{1}8  &  $6.0\times 10^{-8\phantom{0}}$  &  $10^{\pm   38\phantom{0}}$ \\
    binary64  &   53  &  11 &  $1.1\times 10^{-16}$  &  $10^{\pm  308}$ \\
    \bottomrule
    \end{tabularx}
  \end{table}

  As is customary, computed quantities wear a hat in our error analysis.
  Many bounds will feature the constants
  \begin{equation*}
    \gamma_n^h = \frac{n\uhigh}{1-n\uhigh},\qquad
    \gamma_n^\ell = \frac{n\ulow}{1-n\ulow},
  \end{equation*}
  where $n$ is a positive integer.

  We denote the spectrum of a square matrix $A$ by $\Lambda(A)$, and its spectral radius by $\rho(A)$.
  We denote by $\norm{\,\cdot\,}$ any consistent matrix norm.
  For $A\in\Cmn$, we will use the Frobenius norm
  \begin{equation*}
    \lVert A  \rVert_F:=\biggl(\sum_{i=1}^{m}\sum_{j=1}^{n} \lvert a_{ij} \rvert^2\biggr)^{1/2},
  \end{equation*}
  and the induced matrix $p$-norms
  \begin{equation}
    \label{eq:def-p-norm}
    \lVert  A \rVert_p := \max_{x \in \Cn}\frac{\lVert  Ax \rVert_p}{\lVert x \rVert_{p}},\qquad
    \lVert x \rVert_p := \biggl(\sum_{i=1}^{n} \lvert x_i \rvert^p\biggr)^{1/p},\quad
    1 \le p < \infty.
  \end{equation}
  For $p = \infty$, taking the limit of the definition~\eqref{eq:def-p-norm} gives
  \begin{equation*}
  \normo{A} := \lim_{p\to\infty}\norm{A}_p = \max_{1\le i \le m } \sum_{j=1}^{n} \lvert a_{ij} \rvert.
  \end{equation*}
  For a nonsingular matrix $A$ and a vector $x$, we define a normwise condition number, which for the $p$-norms has the form
  \begin{equation*}
    \kappa_{p}(A) := \lVert A \rVert_{p} \cdot \lVert A^{-1} \rVert_{p},
  \end{equation*}
  and componentwise condition numbers~\cite[sect.~7.2]{high:ASNA2}
  \begin{equation*}
    \cond(A) := \normo{\abs{A^{-1}}\abs A}, \quad
    \cond(A,x) := \frac{\normo{\abs{A^{-1}}\abs A\abs x}}{\normo x},
  \end{equation*}
  where the absolute value of a matrix is to be understood componentwise.
  Similarly, inequalities between vectors or matrices are to be interpreted componentwise.

  \subsection{Error analysis}

  The perturbation analysis by Higham~\cite{high93y} shows
  that, for an approximate solution to~\eqref{eq:sylv}, the (normwise) backward
  error~\cite[eq.~(16.10)]{high:ASNA2}
  can be arbitrarily larger than the (normwise) relative residual, which is
  defined by
  \begin{equation}\label{eq:sylv-norm-res}
  	\frac{\norm{AX + XB - C}}
  	{\norm{C} + \norm{X}\bigl(\norm{A} + \norm{B}\bigr)}.
  \end{equation}
  This is in stark contrast with the case of linear systems, for which
  the two quantities are one and the same~\cite[Chap.~7]{high:ASNA2}.
  For two matrices $A \in \Cmm$ and $B \in \Cnn$, the separation is defined as
  \begin{equation}
    \label{eq:sep}
    \sep (A,-B)  := \min_{\norm{X} \neq 0}\frac{\norm{AX  + XB}}{\norm{X}}.
  \end{equation}
  Varah~\cite{vara79} shows that by taking the Frobenius norm
  in~\eqref{eq:sep} we have
  \begin{equation*}
    \sep_F (A,-B)
    = \norm*{(I_n \otimes A + B^T \otimes  I_m)^{-1}}_2^{-1}
    = \sigma_{\min} (I_n \otimes  A + B^T \otimes  I_m),
  \end{equation*}
  which implies that
  \begin{equation*}
    \norm{X}_F \leq \frac{\norm{C}_F}{\sep_F(A,-B)}.
  \end{equation*}

  In our mixed-precision algorithms, the Schur decomposition of $A$ and $B$ are computed by using the QR algorithm~\mbox{\cite[p.~391]{gova13}} in the low precision $\ulow$. We have
  \begin{equation}
  \label{eq:schur-lp}
      A \approx \wUA \wTA \wUA^*,\qquad B \approx \wUB \wTB \wUB^*,
  \end{equation}
  where the computed upper quasi-triangular factors $\wTA \in \Cmm$ and $\wTB \in \Cnn$ satisfy
  \begin{equation}
    \label{eq:bound-err-quasi-triangular-factor}
    \begin{aligned}
      U^* (A+\DA) U = \wTA,\qquad&\normt{\DA} \approx \ulow\normt{A},\\
      V^* (B+\DB) V = \wTB,\qquad&\normt{\DB} \approx \ulow\normt{B},
    \end{aligned}
  \end{equation}
  for some unitary matrices $U \in \Cmm$ and $V \in \Cnn$, and the computed unitary factors $\wUA \in \Cmm$ and $\wUB \in \Cnn$ satisfy
  \begin{equation*}
    \normt*{\wUA^*\wUA - I_m} \approx \ulow,
    \qquad
    \normt*
    {\wUB^*\wUB - I_n} \approx \ulow.
  \end{equation*}

  \subsection{Numerical algorithms}\label{sec:numer-algorithms}

  When working with only one precision, a well-conditioned Sylvester equation can
  be solved accurately and efficiently using the algorithm of Bartels and
  Stewart~\cite{bast72}, of which several implementations are available.
  The algorithm simplifies in the case of the Lyapunov equation~\eqref{eq:lyap}, where the coefficient matrices are adjoints of each other.

  We now recall the Bartels--Stewart algorithm and a cheaper variant that can be derived if the coefficients $A$ and $B$ are Hermitian.

  \subsubsection{The Bartels--Stewart algorithm}
  \label{sec:BS-algorithm}
  The algorithm requires three steps.
  First, we compute the Schur decompositions $A =: \UA \TA \UA^*$ and
  $B =: \UB \TB \UB^*$, where $\TA$ and $\TB$ are upper quasi-triangular
  matrices and $\UA$ and $\UB$ are unitary matrices.
  Next, we multiply~\eqref{eq:sylv} by $\UA^*$ on the left and by $\UB$ on the right to obtain
\begin{equation}
    \label{eq:sylv-tri}
    \TA Y + Y \TB = \Ctri,\qquad
    \Ctri := \UA^* C \UB,
  \end{equation}
  where the coefficients $\TA$ and $\TB$ are upper quasi-triangular.
  Finally, we solve~\eqref{eq:sylv-tri} with a recurrence.
  The result is obtained by observing that if $Y$ satisfies~\eqref{eq:sylv-tri}, then $X:=\UA Y \UB^*$ satisfies~\eqref{eq:sylv}.

  Computing the Schur decomposition of $A$ and $B$ requires approximately $25m^3$ and $25n^3$ flops~\cite[p.~391]{gova13}, respectively, $2mn(m+n)$ flops are needed to compute the updated right-hand side $\Ctri$, solving the Sylvester equation with triangular coefficients requires $mn(m+n)$ flops~\cite[p.~398]{gova13}, and recovering the solution involves two matrix multiplications for $2mn(m+n)$ additional flops.
  Overall, the algorithm requires $25(m^3 + n^3) + 5mn(m+n)$ flops.

  For the Lyapunov equation~\eqref{eq:lyap}, the algorithm can be simplified, as only one Schur decomposition is needed.
  As $m=n$ in this case, the Lyapunov equation can be solved with only
  $35n^3$ flops if the Bartels--Stewart algorithm is used.

  \newcommand{\diagA}{\ensuremath{D_{A}}}
  \newcommand{\diagB}{\ensuremath{D_{B}}}
  \subsubsection{Hermitian matrix coefficients} If the coefficients $A$ and
  $B$ in \eqref{eq:sylv} are both Hermitian, the Bartels--Stewart
  algorithm can be simplified---and its computational cost can be
  significantly reduced---by exploiting the fact that the quasi-triangular
  Schur factor of a Hermitian matrix is diagonal.
  In fact, once the decompositions $A=:\UA \diagA \UA^{*}$ and $B=:\UB
  \diagB \UB^{*}$ are computed, solving the matrix equation
  \begin{equation*}
    \label{eq:sylv-pert-diag}
    \diagA Y + Y\diagB = \Ctri,\qquad
    \Ctri := \UA^* C \UB,
  \end{equation*}
  amounts to applying the formula
  \begin{equation*}
    Y_{ij} = \frac{\Ctri_{ij}}{(\diagA)_{ii} + (\diagB)_{jj}},
  \end{equation*}
  and then retrieving the solution as
  \begin{equation}
    \label{eq:herm-final-step}
    X = \UA Y \UB^{*}.
  \end{equation}

  Computing the eigendecomposition of an $n \times n$ matrix requires
  $9n^3$
  flops, determining $Y$ requires only $2n^2$ flops,
  and
  computing the new right-hand side $\wt C$ and computing $X$ as
  in~\eqref{eq:herm-final-step} have the same
  cost as in the non-Hermitian case.
  Therefore, the algorithm for Hermitian matrices asymptotically requires
  $26n^3$ flops---the Bartels--Stewart method would need about $60n^3$ flops when $m=n$.

  \section{Iterative refinement for perturbed quasi-triangular equations}
  \label{sec:pert-sylv-algs}

  \newcommand{\correction}{\ensuremath{\delta x}}

  Consider the perturbed quasi-triangular Sylvester equation
  \begin{equation}
    \label{eq:sylv-pert-tri}
    (\TA+\DTA)Y + Y(\TB+\DTB) = \Ctri,
  \end{equation}
  where $\TA \in \Cmm$ and $\TB \in \Cnn$ are upper quasi-triangular, and
  $\DTA \in \Cmm$ and $\DTB \in \Cnn$ are \emph{unstructured} matrices.
  The perturbed Sylvester operator can be split into a dominant term $\mathcal{D}: \Cmn\mapsto \Cmn$ and a perturbation term $\mathcal{P}: \Cmn\mapsto \Cmn$ such that
  \begin{equation*}
      \mathcal{D}(Y) + \mathcal{P}(Y) = \wt C,\qquad
      \mathcal{D}(Y) :=  \TA Y + Y\TB,\qquad
      \mathcal{P}(Y) :=  \DTA Y + Y \DTB.
  \end{equation*}
  Thus, we obtain the fixed‐point iteration
  \begin{equation*}
      \mathcal{D}(Y_{i+1}) = \Ctri - \mathcal{P}(Y_i),
  \end{equation*}
  where at each step we solve
  \begin{equation*}
    \TA Y_{i+1} + Y_{i+1} \TB = \Ctri - (\DTA Y_{i} + Y_{i} \DTB).
  \end{equation*}
  Equivalently, by setting $Y_{i+1} := Y_i + D_i$ and rearranging the terms, we can solve the triangular Sylvester equation
  \begin{equation}
    \label{eq:matr-form-iter}
    \TA D_i + D_i \TB = \Ctri - (\TA+\DTA)Y_{i} - Y_{i} (\TB+\DTB)
  \end{equation}
  for the correction $D_i$.
  Using~\eqref{eq:matr-form-iter} repeatedly yields a classical fixed‑point iteration with regular splitting, which represents the foundation of the Alternating Direction Implicit (ADI) method~\cite{luwa91}, \cite{pera55}.

  \newcommand{\funfont}[1]{\textsc{#1}}
  \newcommand{\funstit}{\ensuremath{\funfont{solve\_pert\_sylv\_tri\_stat}}}
  \newcommand{\funref}{\ensuremath{\funfont{solve\_pert\_sylv\_tri}}}

  \begin{algorithm2e}[t]
    \caption{Stationary iteration-like method for the solution
      of~\eqref{eq:sylv-pert-tri}.}
    \label{alg:pert-sylv-tri-stationary}
    \KwIn{Known matrices in~\eqref{eq:sylv-pert-tri},
      initial guess $Y_0$, and target tolerance $\varepsilon>0$.}
    \KwOut{$Y \in \Cmn$ such that $(\TA + \DTA) Y + Y (\TB + \DTB) \approx \Ctri$.}
    \Function{\funstit\textup{(}$\TA$, $\DTA$, $\TB$, $\DTB$, $\Ctri$, $Y_0$\textup{)}}{
      $i \gets 0$\;
      \Repeat{$\norm{D_{i-1}} / \norm{Y_{i}} \le \varepsilon$}{
        Find $D_i$ such that $\TA D_{i} + D_{i} \TB = \Ctri - (\TA + \DTA)
        Y_i-Y_i(\TB + \DTB)$.\;\label{ln:sylv-eq1}
        $Y_{i+1} \gets Y_i + D_i$\;\label{ln:sylv-update1}
        $i \gets i+1$\;
      }
      \Return $Y\gets Y_{i}$\;
    }
  \end{algorithm2e}

  The pseudocode of a possible implementation is given in
  \cref{alg:pert-sylv-tri-stationary}, where the matrix equation on \cref{ln:sylv-eq1} can be solved with the
  Bartels--Stewart substitution algorithm.
  The function \funstit{} will be a building block of our mixed-precision algorithms in later sections.
  For this purpose, we will assume that all computations are performed using the highest precision among the input arguments.

An alternative iterative refinement scheme for the generalized Sylvester
equation is proposed in~\cite[sect.~6.4]{kohl21}, where the equation is
formulated differently---not in a perturbed form. Instead of solving a
quasi-triangular equation, that algorithm directly recovers the solution to the
full equation, requiring additional matrix multiplications with the unitary
factors of the generalized Schur decomposition in each iteration. This approach
does not seem applicable in a mixed-precision environment, as it focuses on the
full equation rather than on a quasi-triangular one.

  The standard convergence theory of fixed-point iteration for linear systems~\cite[sect.~4.2]{saad03} shows that
  \begin{equation}
    \label{eq:norm-2-fro}
    \norm{\DTA}_2 + \norm{\DTB}_2 \le \norm{\DTA}_F + \norm{\DTB}_F <
    \sep_F(\TA,-\TB)
  \end{equation}
  is a sufficient condition for the convergence of \cref{alg:pert-sylv-tri-stationary}.

  If $\DTA$ and $\DTB$ are upper (quasi) triangular, or $\TA$ and $\TB$ commute with $\DTA$ and $\DTB$, respectively, then one can show that this result holds for any \mbox{operator norm}.

\subsection{Computational cost and storage requirement}
  We now assess the computational cost of \cref{alg:pert-sylv-tri-stationary}.
  Each step requires $2mn(m+n)$ flops to compute the right-hand side of the Sylvester equation on \cref{ln:sylv-eq1} and $mn(m+n)$ flops to solve it.
  Overall, \cref{alg:pert-sylv-tri-stationary} requires $3kmn(m+n)$ flops, where $k$ is the total number of iterations required to achieve convergence.
  For Lyapunov equations, the computational cost of the algorithm can be reduced to $6kn^3$~flops.

  Finally, we comment on the storage requirements of \cref{alg:sylv-mp} in terms of floating-point values (flovals), assuming that the computation is performed using the relevant BLAS~\cite{bddd02} and LAPACK~\cite{abbb99} routines.
  Assuming that $\DTA$ and $\DTB$ can be overwritten, storing $\TA + \DTA$ and $\TB + \DTB$ does not require any additional memory.
  Computing the right-hand side of the Sylvester equation requires a temporary $m \times n$ matrix and two applications of \verb|xGEMM|.
  The triangular Sylvester equation can be solved using \verb|xTRSYL|, which does not require any additional memory, or \verb|xTRSYL3|, which requires an amount of additional memory depending on the architecture where the programme is run.
  Both functions overwrite the right-hand side of the equation with the solution.
  Therefore, at a minimum \funstit{} requires $mn$ flovals of additional storage.

    \subsection{Error analysis and attainable residual}
  If the sufficient condition~\eqref{eq:norm-2-fro} is satisfied, then \cref{alg:pert-sylv-tri-stationary} converges in exact arithmetic.
   What can be said about the convergence of this algorithm in floating-point arithmetic?
  To answer this question, we leverage the connection between~\cref{alg:pert-sylv-tri-stationary} and a variant of fixed-precision iterative refinement for linear systems---by analysing the propagation of errors in the latter, we investigate the attainable relative residual of~\cref{alg:pert-sylv-tri-stationary}.
This is possible because, in our mixed-precision setting, the magnitude of the entries of $\DTA$ and $\DTB$ depends on the magnitude of $\ulow$, the unit roundoff of the low precision.

  \begin{algorithm2e}[t]
  	\caption{\mbox{Fixed-precision iterative refinement variant to solve $Mx=b$.}}
  	\label{alg:ir-linear-system}
  	\KwIn{$M$, $\DM\in\Css$, $b\in\Cs$, initial guess $x_0\in\Cs$, target tolerance $\varepsilon>0$.}
  	\KwOut{Approximate solution $\widehat{x}$ to $Mx = b$.}
  	$i \gets 0$\;
  	\Repeat{$\norm{d_{i-1}} / \norm{x_{i}} \le \varepsilon$}{
  		Compute $r_i = b - Mx_i$.\; \label{ln:compt-res}
  		Solve $(M-\DM)d_i = r_i$. \;\label{ln:linear-solve}
  		$x_{i+1} = x_i + d_i$ \; \label{ln:sol-update}
  		$i \gets i+1$\;
  	}
  	\Return $\widehat{x}=x_{i}$\;
  \end{algorithm2e}

  To understand the convergence of~\cref{alg:pert-sylv-tri-stationary} in finite-precision arithmetic, it is sufficient to consider the equivalent iterative refinement process for the perturbed linear system
  \begin{equation}\label{eq:sylv-per-tri-linear}
      (M_t + \DM_t)\vvec(Y) = \vvec(\Ctri),
  \end{equation}
  where
  $$
    M_t := I_n\otimes \TA+\TB^T\otimes I_m, \quad
    \DM_t := I_n\otimes \DTA+\DTB^T\otimes I_m.
  $$
  To facilitate the analysis, we recast this fixed-precision iterative refinement scheme as an algorithm for solving the (unperturbed) linear system $Mx=b$.
  The pseudocode of this method is given in~\cref{alg:ir-linear-system}.
  Our analysis does not exploit the block-triangular structure of $M$, because all bounds we use in the derivation are satisfied by general matrices. We note that this generic treatment on the special structure could lead to a potentially large overestimation on the errors arising in the actual Sylvester setting.
  The difference between \cref{alg:ir-linear-system} and traditional fixed-precision iterative refinement~\cite[p.~232]{high:ASNA2} is on~\cref{ln:linear-solve}, where our variant solves a perturbed linear system for the update~$d_i$.

  Since \cref{alg:pert-sylv-tri-stationary} is only to be used in precision $\uhigh$,
  we assume that the working precision of~\cref{alg:ir-linear-system} is $\uhigh$ throughout.
  To investigate the behavior of the forward and backward errors of the solution computed by~\cref{alg:ir-linear-system}, we assume that the computed solution $\whd_i$ to $(M - \DM)\whd_i  = \whr_i$ on~\cref{ln:linear-solve} satisfies
  \begin{align}
  	\whd_i & = (I_{s}+\uhigh G_i)d_i, \qquad \uhigh \normo{G_i}<1, \label{eq:assum-di} \\
  	\normo{\whr_i - (M-\DM) \whd_i} & \le \uhigh
  	(c_1\normo{M-\DM}\normo{\whd_i} + c_2 \normo{\whr_i}).
  	\label{eq:assum-res}
  \end{align}
  Here, $G_i$, $c_1$, and $c_2$ depend on $s$, $M-\DM$, $\whr_i$,
  and $\uhigh$.
Equations~\eqref{eq:assum-di} and~\eqref{eq:assum-res}, which are similar to~\cite[eqs.~(2.3) and (2.4)]{cahi18}, constrain the relative forward and backward errors, respectively: the first requires that the forward error be bounded by a multiple of $\uhigh$ strictly less than $1$; the second requires that the backward error~\cite[Thm.~7.1]{high:ASNA2} be of order at most $\max(c_1,c_2)\uhigh$.
  Also, we assume that the perturbation $\DM$ in~\cref{alg:ir-linear-system} satisfies
  \begin{equation}\label{eq:assum-pertb}
  	\normo{\DM} \le c_3 \ulow \normo{M}, \quad
  	\rho(M^{-1}\DM)<1, \quad
  	0\le c_3\equiv c_3(M, \DM)\le s,
  \end{equation}
  where $s$ is the order of $M$.
  This assumption on the unstructured perturbation is a consequence of the equivalence between~\eqref{eq:sylv-pert-tri} and~\eqref{eq:sylv-per-tri-linear}, together with~\eqref{eq:bound-err-quasi-triangular-factor}, which bounds the perturbation if the triangular Schur factors of $A$ and $B$ are computed in low precision.
  As the bound~\eqref{eq:bound-err-quasi-triangular-factor} is given in the $2$-norm, we have introduced the constant $c_3$ to account for the shift of norm.
  A sufficient condition for the spectral radius bound in~\eqref{eq:assum-pertb}, which is to guarantee the convergence of the fixed-point iteration for solving~\eqref{eq:sylv-per-tri-linear}, is $c_3 \ulow \kappa_{\infty}(M)<1$. This bound will play a crucial role in our characterization of the forward and backward errors (see \cref{thm:forward-converge-accuracy} and \cref{thm:backward-converge-accuracy}), and it essentially requires that $M$ is not too ill conditioned with respect to precision~$\ulow$.

  \newcommand{\DMb}{\Delta \overline{M}}

  \subsection{Normwise forward error analysis}
  We begin by analyzing the behavior of the forward error of $\widehat{x}_i$, the approximated solution produced by~\cref{alg:ir-linear-system} at the $i$th iteration.

  For the computation of $r_i$ on \cref{ln:compt-res} of~\cref{alg:ir-linear-system}, one can show that~\cite[sect.~12.1]{high:ASNA2}
  \begin{equation}\label{eq:forward-delta-res}
    \widehat{r}_i = b-M\widehat{x}_i + \Delta r_i,\quad \lvert \Delta r_i \rvert \le
    \gamma_{s+1}^h\big(\lvert b \rvert + \lvert M \rvert\lvert \widehat{x}_i \rvert\big) \le
    \gamma_{s+1}^h \big(\lvert M \rvert \lvert x-\widehat{x}_i \rvert + 2\lvert M \rvert\lvert x \rvert\big).
  \end{equation}
  We note that
  \begin{equation}\label{eq:steps-41}
  (M-\DM)^{-1} \whr_i
  		\approx (I + M^{-1}\DM)M^{-1} \whr_i
  		= (I + M^{-1}\DM) (x-\widehat{x}_i + M^{-1}\Delta r_i),
  \end{equation}
  where we have used the approximation
  \begin{equation}\label{eq:pert-inv-expan}
  	(M-\Delta M)^{-1} \approx (I+M^{-1}\Delta M) M^{-1}.
  \end{equation}
  Equation~\eqref{eq:pert-inv-expan} is obtained by ignoring higher-order terms in the Neumann expansion of $(I-M^{-1}\DM)^{-1}$, which is guaranteed to converge in view of~\eqref{eq:assum-pertb}.
  It follows, by using~\eqref{eq:assum-di}, that the solver on \cref{ln:linear-solve} satisfies
    \begin{equation}\label{eq:steps-42}
  \whd_i - (M-\DM)^{-1} \whr_i
   = \uhigh G_i(I + M^{-1}\DM) (x-\widehat{x}_i + M^{-1}\Delta r_i).
  \end{equation}
   By substituting~\eqref{eq:forward-delta-res} into~\eqref{eq:steps-42}, we obtain
   \begin{align}\label{eq:forward-bound}
    \lvert \whd_i - (M-\DM)^{-1} \whr_i \rvert \le&
    \uhigh \lvert G_i \rvert\lvert I + M^{-1}\DM \rvert \big(\lvert \wh{x}_i-x \rvert +
    \gamma_{s+1}^h \lvert M^{-1} \rvert\lvert M \rvert(\lvert x-\wh{x}_i \rvert+2\lvert x \rvert) \big) \nonumber\\
    \le& \uhigh \lvert G_i \rvert\lvert I + M^{-1}\DM \rvert(I + \gamma_{s+1}^h \lvert M^{-1} \rvert\lvert M \rvert)\lvert \wh{x}_i-x \rvert\nonumber\\
    &+ 2\uhigh \gamma_{s+1}^h\lvert G_i \rvert\lvert I + M^{-1}\DM \rvert \lvert M^{-1} \rvert\lvert M \rvert\lvert x \rvert.
   \end{align}
 Following the modified standard model of floating-point arithmetic~\cite[eq.~(2.5)]{high:ASNA2}, we see that the solution update on \cref{ln:sol-update} satisfies
  \begin{equation}\label{eq:forward-delta-xi}
    \widehat{x}_{i+1} =
    \widehat{x}_i + \widehat{d}_i + \Delta {x}_i, \quad
    \lvert \Delta {x}_i \rvert
    \le \uhigh\lvert \widehat{x}_{i+1} \rvert.
  \end{equation}
  By substituting~\eqref{eq:steps-41} into~\eqref{eq:forward-delta-xi}, we arrive at
  \begin{align}\label{eq:forward-xi+1}
    \widehat{x}_{i+1}
    & = \widehat{x}_i + (M-\DM)^{-1} \whr_i +
   \big(\widehat{d}_i - (M-\DM)^{-1} \whr_i + \Delta {x}_i\big)\nonumber\\
   & \approx x + M^{-1}\Delta r_i + M^{-1}\Delta M(x-\widehat{x}_i + M^{-1}\Delta r_i) + \big(\widehat{d}_i - (M-\DM)^{-1} \whr_i + \Delta {x}_i\big),
  \end{align}
  and, therefore, from~\eqref{eq:forward-bound} and \eqref{eq:forward-delta-xi}, we conclude that
  \begin{align*}
    \lvert \wh{x}_{i+1}-x \rvert &\lesssim \lvert M^{-1} \rvert\lvert I + \DM M^{-1} \rvert\cdot \gamma_{s+1}^h(\lvert M \rvert \lvert x-\widehat{x}_i \rvert + 2 \lvert M \rvert\lvert x \rvert) + \lvert M^{-1}\DM \rvert\lvert \wh{x}_{i}-x \rvert  \\
    &\phantom{\le 3} + \uhigh \lvert G_i \rvert\lvert I + M^{-1}\DM \rvert(I + \gamma_{s+1}^h \lvert M^{-1} \rvert\lvert M \rvert)\lvert \wh{x}_i-x \rvert \\
    &\phantom{\le 3} + 2\uhigh \gamma_{s+1}^h\lvert G_i \rvert\lvert I + M^{-1}\DM \rvert \lvert M^{-1} \rvert\lvert M \rvert\lvert x \rvert + \uhigh\lvert \widehat{x}_{i+1} \rvert  \\
    &=: F_i \abs{\widehat{x}_{i}-x} + f_i,
  \end{align*}
  where
  \begin{align*}
    F_i &= \gamma_{s+1}^h \lvert M^{-1} \rvert\lvert I + \DM M^{-1} \rvert \lvert M \rvert + \lvert M^{-1}\DM \rvert
    + \uhigh \lvert G_i \rvert\lvert I + M^{-1}\DM \rvert(I + \gamma_{s+1}^h \lvert M^{-1} \rvert\lvert M \rvert),\\
    f_i &= 2\gamma_{s+1}^h \bigl( \lvert M^{-1} \rvert\lvert I + \DM M^{-1} \rvert + \uhigh \lvert G_i \rvert\lvert I + M^{-1}\DM \rvert \lvert M^{-1} \rvert \bigr)\lvert M \rvert \lvert x \rvert + \uhigh\lvert \widehat{x}_{i+1} \rvert.
  \end{align*}
  To shorten the notation in the equations in this and the following section, we define
  \begin{equation}
    \label{eq:psi-M}
    \psi_{M}:= 1 +  c_3\ulow\kappa_{\infty}(M)
  \end{equation}
  By using the bound $\gamma_{s+1}^h \le (s+1)\uhigh$ and~\eqref{eq:assum-pertb}, it is straightforward to show that
  \begin{align*}
    \normo{F_i} & \le (\uhigh \normo{G_i}+1)
    \psi_{M}( 1 + (s+1)\uhigh\cond(M)) -1 \\
    & \less 2\psi_{M}( 1 + (s+1)\uhigh\cond(M)) -1
  \end{align*}
  and that
  \begin{align*}
    \normo{f_i} & \le 2(s+1)\uhigh \bigl( \psi_{M} +
    \uhigh \normo{G_i} \psi_{M} \bigr) \normo{\lvert M^{-1} \rvert \lvert M \rvert \lvert x \rvert} + \uhigh\normo{\widehat{x}_{i+1}} \\
    & \less 4(s+1)\uhigh \psi_{M}
    \normo{\lvert M^{-1} \rvert \lvert M \rvert \lvert x \rvert}  + \uhigh\normo{\widehat{x}_{i+1}}.
  \end{align*}
  We summarize our findings in the following result.

  \begin{theorem}\label{thm:forward-converge-accuracy}
    Let~\cref{alg:ir-linear-system} be applied to a linear system $Mx = b$ in precision $\uhigh$, where $M\in\mathbb{C}^{s\times s}$ is nonsingular, and assume that the solver used on \cref{ln:linear-solve} satisfies~\eqref{eq:assum-di} and~\eqref{eq:assum-pertb}.
    If
    $$
    \phi_i := (1+ \uhigh \normo{G_i})
    \bigl(1 +  c_3\ulow\kappa_{\infty}(M)\bigr)\bigl( 1 + (s+1)\uhigh\cond(M)\bigr)
    $$
    is sufficiently smaller than $2$ \textup{(}with the precise bound depending on $\normo{f_i}$\textup{)},
    then, at iteration $i$, the normwise forward error is reduced by a factor approximately $\phi_i$, and this behavior persists until an iterate $\widehat x$ is produced for which
    $$
    \frac{\normo{x-\widehat x}}{\normo{x}} \lesssim
    4(s+1)\uhigh \bigl(1+c_3 \ulow \kappa_{\infty}(M)\bigr)
    \cond(M,x)  + \uhigh.
    $$
  \end{theorem}

  The $\phi_i$ in \cref{thm:forward-converge-accuracy} is a product whose terms depend on both $\uhigh$ and $\ulow$, and the most vulnerable factor in the rate of convergence is $\psi_{M}$ in~\eqref{eq:psi-M}: a necessary condition for $\phi_i<2$ to hold is that $c_3\ulow\kappa_{\infty}(M)<1$.
  The limiting accuracy of~\cref{alg:ir-linear-system} is a factor $\psi_{M}$ worse than for the standard fixed-precision iterative refinement (cf.~\cite[Cor.~3.3]{cahi18}): this is a consequence of the noise $\DM$, which has fixed magnitude at each step.

  We repeated this analysis following the approach of Higham~\cite{high97i},~\cite[sect.~12.1]{high:ASNA2} and, with similar assumptions, we obtained essentially the same bound for the limiting accuracy of the forward error.

  \subsection{Normwise backward error analysis}
  Now we turn our focus to the behavior of the normwise backward error, that is, the normwise residual.

  Multiplying~\eqref{eq:forward-xi+1} by $M$, and using~\eqref{eq:steps-41} and then~\eqref{eq:forward-delta-res}, gives
\begin{align}\label{eq:backward-key}
    M\wh{x}_{i+1}-b &\approx  \Delta r_i + \Delta M(x-\widehat{x}_i + M^{-1}\Delta r_i) + M\widehat{d}_i - M(M-\DM)^{-1} \whr_i + M\Delta {x}_i \nonumber\\
    &= \Delta r_i + (M-\DM)\widehat{d}_i - \whr_i +  M\Delta {x}_i
    + \DM \whd_i.
  \end{align}
  Let $h_i :=(M-\DM)\whd_i - \whr_i$. If we take the norm of $\whd_i = (M-\DM)^{-1}(h_i + \whr_i)$, using assumptions~\eqref{eq:assum-res} and \eqref{eq:assum-pertb} and the quantity defined in~\eqref{eq:psi-M}, we obtain
  \begin{align*}
    \normo{h_i} & \le \uhigh
    \big(c_1\normo{M-\DM}\normo{\whd_i} + c_2 \normo{\whr_i}\big) \\
    & \le \uhigh
    \big(c_1(1+c_3 \ulow) \psi_{M}
    \kappa_{\infty}(M)(\normo{h_i} + \normo{\whr_i}) + c_2 \normo{\whr_i}
    \big),
  \end{align*}
  where we used $\normo{(M-\DM)^{-1}} \le \psi_{M}\normo{M^{-1}}$, which is a consequence of~\eqref{eq:assum-pertb} and~\eqref{eq:pert-inv-expan}.
  Now, assuming that
  \begin{equation*}
    c_4\kappa_{\infty}(M)\uhigh <1, \quad
    c_4 \equiv c_4(M,\ulow,c_1,c_3) :=
    c_1(1+c_3 \ulow)\psi_{M},
  \end{equation*}
  we have
  \begin{equation}\label{eq:backward-hi-norm}
    \normo{h_i} \le \uhigh \frac{c_4\kappa_{\infty}(M) + c_2}
    {1-c_4\kappa_{\infty}(M)\uhigh} \normo{\whr_i}.
  \end{equation}
  If we write $ d_i =(M-\DM)^{-1}\whr_i$, then by assumption~\eqref{eq:assum-di} we get
  \begin{equation}\label{eq:backward-di-hat}
    \normo{\whd_i} \le (1+\uhigh\normo{G_i})
      \psi_{M}\normo{M^{-1}} \normo{\whr_i}.
  \end{equation}
  Therefore, from~\eqref{eq:backward-key}, using~\eqref{eq:forward-delta-xi}, \eqref{eq:backward-hi-norm}, \eqref{eq:backward-di-hat},
  and two invocations of~\eqref{eq:forward-delta-res}, we obtain
  \begin{align*}
    \normo{b-M\wh{x}_{i+1}} \lesssim
    & \gamma_{s+1}^h\big(\normo{b} + \normo{M}\normo{\widehat{x}_i}\big)
      + \uhigh \frac{c_4\kappa_{\infty}(M) + c_2}{1-c_4\kappa_{\infty}(M)\uhigh} \normo{\whr_i} \\
    &+ \uhigh\normo{M}\normo{\wh{x}_{i+1}} +
    c_3 \ulow \kappa_{\infty}(M)
    (1+\uhigh\normo{G_i}) \psi_{M} \normo{\whr_i},
  \end{align*}
  and we finally can conclude that
  \begin{equation*}
     \normo{b-M\wh{x}_{i+1}}  \lesssim
    \psi_i \normo{b-M\wh{x}_{i}}
     +  \xi_i,
  \end{equation*}
  where we have defined
  \begin{align}\label{eq:psi_i}
  	\psi_i &:= \uhigh \frac{c_4\kappa_{\infty}(M) + c_2}{1-c_4\kappa_{\infty}(M)\uhigh} +c_3 \ulow \kappa_{\infty}(M)
  	(1+\uhigh\normo{G_i})\psi_{M}, \\ \nonumber
    \xi_i &:= (s+1)( 1 + \psi_i)\uhigh
    \big(\normo{b} +
    \normo{M}\normo{\widehat{x}_i}\big)
     + \uhigh\normo{M}\normo{\wh{x}_{i+1}}
  \end{align}
  and have used the bound $\gamma_{s+1}^h \le (s+1)\uhigh$.

  We now conclude the analysis above and state a result on the
  behavior of the normwise backward error and the limiting residual of~\cref{alg:ir-linear-system}.

  \begin{theorem}\label{thm:backward-converge-accuracy}
    Let~\cref{alg:ir-linear-system} be applied to a linear system $Mx = b$ in precision $\uhigh$, where $M\in\mathbb{C}^{s\times s}$ is a nonsingular matrix satisfying $c_4\kappa_{\infty}(M)\uhigh <1$ with
    $c_4 = c_1(1+c_3 \ulow)\bigl(1+c_3 \ulow \kappa_{\infty}(M)\bigr)$.
    Assume that the solver used on \cref{ln:linear-solve} satisfies~\eqref{eq:assum-di}--\eqref{eq:assum-pertb}.
    If $\psi_i$ in~\eqref{eq:psi_i}
    is sufficiently smaller than $1$ \textup{(}with the precise bound depending on $\xi_i$\textup{)}, then, at iteration $i$, the normwise residual is reduced by a factor approximately $\psi_i$, and this behavior persists until an iterate $\widehat x$ is produced for which
    $$
    {\normo{b-M\widehat x}} \lesssim
    2(s+1)\uhigh\bigl(\normo{b} +
    \normo{M}\normo{\widehat{x}}\bigr)  + \uhigh\normo{M}\normo{\wh{x}}.
    $$
  \end{theorem}

  Again, $c_3\ulow\kappa_{\infty}(M)<1$ is a necessary condition for $\psi_i<1$ to hold, and the convergence condition for the backward error is more stringent than that for the forward error.
  This is not surprising, since backward stability is a stronger property than forward stability.
  Under the conditions of \cref{thm:backward-converge-accuracy}, we will eventually produce an $\wh x$ such that
  \begin{equation*}
    {\normo{b-M\widehat x}} \lesssim 2(s+1)\uhigh
    \bigl(\normo{b} + \normo{M}\normo{\widehat{x}}\bigr),
  \end{equation*}
  which implies that $\wh{x}$ is a backward stable solution to precision $\uhigh$~\cite[Thm.~7.1]{high:ASNA2}.

 \Cref{tab:ir-linear-system} summarizes the results in \cref{thm:forward-converge-accuracy,thm:backward-converge-accuracy} and provides the limiting accuracy of~\cref{alg:ir-linear-system} under different choices of precisions and on problems with varying condition numbers. The table is comparable with~\cite[Tab.~7.1]{cahi18}
  (only in the case of $u_s=u_f=\ulow$ and $u=u_r=\uhigh$ therein),
  because, with our assumptions~\eqref{eq:assum-di}--\eqref{eq:assum-pertb},~\cref{alg:ir-linear-system} can be thought of as solving the linear system $Mx=b$ on \cref{ln:linear-solve} accurately in precision $\ulow$.

  \begin{table}[t]
    \centering
    \caption{Different choices of floating-point arithmetics for~\cref{alg:ir-linear-system}.
      The third column shows an approximate bound on $\kappa_\infty (M)$ that must hold for the analysis to guarantee convergence with limiting backward or forward errors of the orders shown in the final two columns.}
    \label{tab:ir-linear-system}
    \begin{tabularx}{\linewidth}{Xcccc}
      \toprule
      &&& \multicolumn{2}{c}{Limiting error} \\
      \cmidrule{4-5}
      \multicolumn{1}{c}{$\ulow$} & $\uhigh$ & $\kappa_\infty (M)$ & Backward & Forward  \\
      \midrule
      bfloat16                     & binary32 & $10^3$ & $2^{-24}$ & $\cond(M,x)\times 2^{-24}$  \\
      binary16/TensorFloat-32      & binary32 & $10^4$ & $2^{-24}$ & $\cond(M,x)\times 2^{-24}$  \\
      \midrule
      bfloat16                     & binary64 & $10^3$ & $2^{-53}$ & $\cond(M,x)\times 2^{-53}$  \\
      binary16/TensorFloat-32      & binary64 & $10^4$ & $2^{-53}$ & $\cond(M,x)\times 2^{-53}$  \\
      binary32                     & binary64 & $10^8$ & $2^{-53}$ & $\cond(M,x)\times 2^{-53}$  \\
      \bottomrule
    \end{tabularx}
  \end{table}

  Because of the equivalence between the residual of the linear system $Mx=b$ and the residual of the Sylvester equation~\eqref{eq:sylv-pert-tri}, \cref{tab:ir-linear-system} implies that, under one of the presented settings for $\ulow$, $\uhigh$, and $\kappa_\infty (M)$, the \textit{relative residual} of the Sylvester equation~\eqref{eq:sylv-pert-tri} is of the order of the unit roundoff of binary32 or binary64 arithmetic.
  Yet, not much can be said about the backward or forward error of the Sylvester equation without further assumptions~\cite[sect.~16.2]{high:ASNA2}.

  \section{Mixed-precision algorithms for the Sylvester equation}
  \label{sec:mixed-precision}

  We can use the iterative refinement scheme \funstit\ in \cref{alg:pert-sylv-tri-stationary} to develop mixed-precision algorithms for solving the Sylvester equation~\eqref{eq:sylv}.
  With the notation in~\eqref{eq:schur-lp}, let
  \begin{equation}
  \label{eq:rw-pert}
  \wDA = A - \wUA \wTA \wUA^*\quad \text{ and }\quad
  \wDB = B - \wUB \wTB \wUB^*
  \end{equation}
  denote the errors in the Schur decompositions of $A$ and $B$ computed in precision \ulow.
  One might be tempted to use the function \funstit\ to solve
  \begin{equation}
    \label{eq:sylv-mp}
    \bigl(\wTA + \wUA^* \wDA \wUA\bigr) Y +
    Y \bigl(\wTB + \wUB^* \wDB \wUB\bigr) =
    \wUA^* C \wUB
  \end{equation}
  in precision $\uhigh$ and then recover the solution of~\eqref{eq:sylv} as $X\gets \wUA Y \wUB^*$.
  This would not work, because $\wUA$ and $\wUB$ are only unitary to precision \ulow.
  To recover the solution to~\eqref{eq:sylv} from a solution to~\eqref{eq:sylv-tri}, $\UA$ and $\UB$ must the unitary to precision $\uhigh$, and this is not the case if the Schur decomposition is only computed in precision $\ulow$.

  In \cref{sec:mixed-precision-orth} and \cref{sec:mixed-precision-inv}, we develop two mixed-precision algorithms to solve~\eqref{eq:sylv}.
  The former is based on re-orthonormalization in high precision, the other is based on explicit inversion of the almost-unitary factors.

  \newcommand{\statvar}{\ensuremath{{\texttt{stationary}}}}
  \subsection{Orthonormalization of the unitary factors in high
    precision}
    \label{sec:mixed-precision-orth}

    \begin{algorithm2e}[t]
    \caption{\mbox{Orthonormalization-based mixed-precision Sylvester solver.}}
    \algorithmfootnote{On \cref{ln:qr-A,ln:qr-B}, all entries along the diagonal of $\RA$ and $\RB$ must be positive.}
    \label{alg:sylv-mp}
    \KwIn{$A\in\Cmm$, $B\in\Cnn$, $C\in\Cmn$, $\ulow$ and $\uhigh$ as in~\eqref{eq:prec-inequality}.}
    \KwOut{$X \in \Cmn$ such that $A X + X B \approx C$.}
    Compute Schur decomposition $A =: \wUA \wTA \wUA^*$ in precision
    $\ulow$.\;\label{ln:schurA}
    Compute Schur decomposition $B =: \wUB \wTB \wUB^*$ in precision
    $\ulow$.\;\label{ln:schurB}
    Compute QR factorization $\wUA =: \QA \RA$ in precision $\uhigh$.\;\label{ln:qr-A}
    Compute QR factorization $\wUB =: \QB \RB$ in precision $\uhigh$.\;\label{ln:qr-B}
    $F \gets \QA^* C \QB$, computed in precision
    $\uhigh$.\;\label{ln:preproc-begin-1}
    $\LA \gets \QA^* A \QA - \wTA$, computed in precision
    $\uhigh$.\;\label{ln:preproc-mid-1}
    $\LB \gets \QB^* B \QB - \wTB$, computed in precision $\uhigh$.\;\label{ln:preproc-end-1}
    Find $Y_0$ such that $\wTA Y_0 + Y_0 \wTB = F$ in precision
    $\ulow$.\;\label{ln:initial-sylv-1}
    $Y \gets$ \funstit($\wTA$, $\LA$, $\wTB$, $\LB$, $F$, $Y_0$) in precision $\uhigh$.\;
    $X \gets \QA Y \QB^*$ in precision
    $\uhigh$.\;\label{ln:recovery}
  \end{algorithm2e}

  Let $\QA \in \Cmm$ and $\QB \in \Cnn$ be the matrices obtained by
  orthonormalizing, in precision $\uhigh$, $\wUA$ and $\wUB$,
  respectively.
  This can be done by using, for example, the modified Gram--Schmidt (MGS) algorithm, so we have~\cite[Thm.~19.13]{high:ASNA2}, for some orthogonal matrices (in exact arithmetic) $Q_1$ and $Q_2$, that
  \begin{equation}
    \begin{aligned}
    \label{eq:qr-orth-diff-bound}
        \EA:&=\QA-Q_1, \quad
    \normt{\EA}\lesssim d_1\uhigh\kappa_2(\wUA) \approx d_1\uhigh,\\
        \EB:&=\QB-Q_2, \quad
    \normt{\EB}\lesssim d_2\uhigh\kappa_2(\wUB) \approx d_2\uhigh,
    \end{aligned}
  \end{equation}
  and
  \begin{equation*}
    \normt{\QA^* \QA - I} \lesssim d_3 \kappa_2(\wUA)  \uhigh
    \approx d_3 \uhigh,  \quad
    \normt{\QA^* \QA - I} \lesssim d_4 \kappa_2(\wUA)  \uhigh
    \approx d_4 \uhigh,
  \end{equation*}
  for some constants $d_1\equiv d_1(m)$, $d_2\equiv d_2(n)$, $d_3\equiv d_3(m)$, and $d_4\equiv d_4(n)$.
  We used the MGS algorithm because it admits simple and explicit bounds on the loss of orthogonality, Alternatives such as Householder QR and Cholesky-based QR could also be of interest~\cite[Chap.~19]{high:ASNA2}; in practice, Householder QR is often more robust, while Cholesky-based methods can better exploit BLAS-3 operations and be more efficient when the input matrix is nearly unitary.
  Note that $\wUA$ and $\wUB$ in~\eqref{eq:qr-orth-diff-bound} are well conditioned in precision $\uhigh$, which can be assumed since they are orthonormal to precision $\ulow$.
  Moreover, by noting that
  $\normt{\wUA-\QA} \le \normt{\wUA-Q} + \normt{Q- \QA} \lesssim \ulow + d_1\uhigh$, we have
  \begin{align}
  \label{eq:eq:qr-orth-prod-diff}
  G_A := \QA^{*}\wUA-I, \quad
  \normt{G_A} &=
  \normt{\QA^{*} (\wUA-\QA +\QA) -I} \nonumber \\
  & \le \normt{\QA^{*} \QA -I} + \normt{\QA^{*}}\normt{\wUA-\QA}
  \lesssim \ulow +(d_1+d_3)\uhigh.
  \end{align}
  Similarly, we can show $\normt{\QB^{*}\wUB-I} \lesssim \ulow +(d_2+d_4)\uhigh$.
  It requires asymptotically $2m^3$ or $2n^3$ flops to compute a matrix $\QA$ or $\QB$. Once $\QA$ and $\QB$ are obtained, we can use \funstit\ in
  \cref{alg:pert-sylv-tri-stationary} to solve in precision
  $\uhigh$ the perturbed equation
  \begin{equation}
    \label{eq:sylv-mp-orth}
    \bigl(\wTA + (\QA^*A\QA-\wTA) \bigr)Y +
    Y \bigl( \wTB + (\QB^*B\QB-\wTB) \bigr) =
    \QA^*C\QB,
  \end{equation}
  and then recover the solution as $\QA Y \QB^*$.
  This approach is summarized in \cref{alg:sylv-mp}.

  To see why \cref{alg:sylv-mp} computes an approximate solution to the Sylvester equation~\eqref{eq:sylv}, one first should note that~\eqref{eq:sylv-mp-orth}, the equation being solved, is mathematically equivalent to
  \begin{equation}
    \label{eq:sylv-mp-orth-equiv}
    (\QA^*A\QA) Y + Y (\QB^*B\QB) = \QA^*C\QB.
  \end{equation}
  On the other hand,
  with $Y\equiv \QA^*X\QB$ (so $X\equiv \QA^{-*}Y\QB^{-1}$), the Sylvester equation
  \eqref{eq:sylv}
  is mathematically equivalent to
  \begin{equation}\label{eq:sylv-rw-pert2}
    (\QA^*A\QA^{-*} )Y + Y (\QB^{-1}B\QB) =   \QA^*C\QB.
  \end{equation}
  From~\eqref{eq:qr-orth-diff-bound} we have
  \begin{align*}
    \normt{Q_2^{-1}-\QB^{-1}} & = \normt{Q_2^{-1}-(Q_2+\EB)^{-1}} =
    \normt{Q_2^{-1}(I-(I+\EB Q_2^*)^{-1})} \\
    & = \normt{I-(I+\EB Q^*)^{-1}} =
    \normt{\EB Q_2^*-(\EB Q_2^*)^2+\cdots\;}
    \approx\normt{\EB}.
  \end{align*}
  and hence
  \begin{equation*}
    \normt{\QB^{*}-\QB^{-1}} = \normt{\QB^{*}-Q_2^* + Q_2^*-\QB^{-1}}
    \le \normt{\QB -Q_2} +\normt{Q_2^{-1}-\QB^{-1}}
    \approx 2\normt{\EB} \le 2d_2 \uhigh.
  \end{equation*}
  In a similar way we can show that
  $\normt{\QA-\QA^{-*}}=\normt{\QA^*-\QA^{-1}}\lesssim 2d_1 \uhigh$.
  Therefore, \eqref{eq:sylv-mp-orth-equiv} and \eqref{eq:sylv-rw-pert2} differ by two perturbations of order $O(\uhigh)$ in the coefficient matrices on the left-hand side.
  Moreover, the transformation we use to recover the solution $X$ only introduces a perturbation of order $O(\uhigh)$ in norm.
  The whole process amounts to solving in precision $\uhigh$ a nearby equation to~\eqref{eq:sylv} with $O(\uhigh)$ perturbations in the coefficient matrices and in the solution itself.

  Recall that a sufficient condition for the convergence of \cref{alg:pert-sylv-tri-stationary} on the quasi-triangular equation~\eqref{eq:sylv-mp-orth} is the one given in~\eqref{eq:norm-2-fro}. We have
  \begin{equation*}
    \norm{\DTA}
    = \norm{\QA^* A \QA - \wTA}
    = \norm{\QA^{*}(\wUA \wTA \wUA^{*} + \wDA) \QA - \wTA},
  \end{equation*}
  and from~\eqref{eq:eq:qr-orth-prod-diff},
  \begin{align*}
    \norm{\DTA}  &= \norm{ (I+G_A)\wTA(I+G_A^*)-\wTA + \QA^{*}\wDA\QA} \\
    &\lesssim \bigl(\norm{G_A}+\norm{G_A^*}\bigr) \norm{\wTA} + \kappa(\QA)\norm{\wDA}
    \lesssim \ulow\big(2(d_1+d_3+1)\norm{\wTA} + \kappa(\QA)\norm{A} \big).
  \end{align*}
   Similarly, one can show that $\norm{\DTB} \lesssim \ulow\big(2(d_2+d_4+1)\norm{\wTB} + \kappa(\QB)\norm{B} \big)$.

  \begin{table}[t]
    \centering
    \caption{Computational cost (in high-precision and low-precision
      flops) of the algorithms in \cref{sec:mixed-precision}.
      The three constants $\alpha := m^3+n^3$, $\beta := mn(m+n)$, and
      $\gamma := n^{3}$ are used to simplify the notation.}
  \label{tab:mp-comp-cost}
    \begin{tabularx}{\linewidth}{p{2.3cm}YY@{\hskip2pt}cYY}
      \toprule
      & \multicolumn{2}{c}{Sylvester} && \multicolumn{2}{c}{Lyapunov}\\
      \cmidrule{2-3} \cmidrule{5-6}
      & \ulow{} flops & \uhigh{} flops && \ulow{} flops & \uhigh{}
      flops \\
      \midrule
      \cref{alg:sylv-mp}
      & $25\myalpha+\mybeta$ & $6\myalpha+(4+3i)\mybeta$
      && $27\mygamma$ & $(14+6i)\mygamma$ \\
      \cref{alg:sylv-mp-inversion}
      & $25\myalpha+\mybeta$ & $4\frac{2}{3}\myalpha+(4+3i)\mybeta$
      && $27\mygamma$ & $(12\frac{2}{3}+6i)\mygamma$ \\
      \bottomrule
    \end{tabularx}
  \end{table}

  \subsubsection{Computational cost and storage requirement}
  \label{sec:mixed-precision-orth-cost}
  We now discuss the cost of \cref{alg:sylv-mp}.
  The asymptotic computational costs of all the algorithms in this
  section are summarised in \cref{tab:mp-comp-cost}.
  Computing the two Schur decompositions on
  \cref{ln:schurA,ln:schurB} requires $25(m^{3}+n^{3})$
  flops in precision $\ulow$.
  The orthogonal factors of the QR factorizations on
  \cref{ln:qr-A,ln:qr-B} can be computed with
  $2(m^{3}+n^{3})$ flops in precision $\uhigh$ by using the modified
  Gram--Schmidt algorithm, as the matrices $\UA$ and $\UB$ are
  orthonormal to precision $\ulow$ and are thus well conditioned in
  precision $\uhigh$.
  The final recovery of the solution, performed on
  \cref{ln:recovery}, requires two additional matrix products, for
  an additional $2mn(m+n)$ flops in precision~$\uhigh$.

  The computation of $\LA$, $\LB$, and $F$ on \cref{ln:preproc-begin-1,ln:preproc-mid-1,ln:preproc-end-1} requires $4(m^3+n^3) +
  2mn(m+n)$ flops in precision $\uhigh$.
  This algorithm only requires that $Y_{0}$ be an approximate solution
  to the low-precision equation, thus we can solve the Sylvester equation in
  low precision, at the cost of $mn(m+n)$ flops in precision $\ulow$.
  The call to \funstit{} requires $3imn(m+n)$ flops in precision
  $\uhigh$, where $i$ is the number of iterations required by the function
  that solves the triangular Sylvester equation.
  Overall, \cref{alg:sylv-mp} requires $25(m^3+n^3) + mn(m+n)$
  low-precision flops and $6(m^3+n^3)+(4+3i)mn(m+n)$ high-precision ones.
  For a Lyapunov equation, only one Schur decomposition and one
  orthonormalization in high precision are necessary, and $\LB = \LA^*$.
  Therefore, \cref{alg:sylv-mp} requires only $27n^3$ low-precision
  flops and $(14+6i)n^3$ high-precision ones in this case.
  In this scenario, the Bartels--Stewart algorithm would require
  $25(m^{3}+n^{3})+5mn(m+n)$ high-precision flops to solve a Sylvester
  equation
  and $35n^3$ high-precision flops to solve a Lyapunov equation.

  We now assess the additional storage needed by \cref{alg:sylv-mp} in terms of floating-point values (flovals) in precisions $\ulow$ and $\uhigh$.
  To compute the Schur decomposition of $A$ and $B$, we need to convert these matrices to precision $\ulow$, compute their Schur decomposition, and convert the resulting matrices back to high precision.
  The storage requirements can be significantly reduced if the two matrices are considered one at a time.
  Converting the larger of $A$ and $B$ to low precision requires $\max(m,n)^{2}$ flovals in precision $\ulow$.
  The Schur decomposition can be computed using \verb|xGEES|, which overwrites its input with the triangular Schur factor and requires additional storage for the unitary factor.
  Therefore, \cref{ln:schurA,ln:schurB} require $\max(m,n)^{2}$ additional flovals in precision $\ulow$, in addition to the work memory required by the routine, which is at least $3\max(m,n)$ flovals.
  Therefore, computing the Schur factors in precision $\ulow$ and converting them to precision $\uhigh$ requires $\max(m,n)^{2} + 3\max(m,n)$ flovals in precision $\ulow$ and $2(m^{2} + n^{2})$ flovals in precision $\uhigh$ to store $\wUA$, $\wTA$, $\wUB$, and $\wTB$.

  To compute the QR factorization of the unitary factors of the Schur decomposition, one can use \verb|xGEQRF|, which overwrites its input and require at least $\max(m,n)$ flovals of work memory in precision $\uhigh$.
  Assuming that the inputs $A$, $B$, and $C$ can be overwritten, applying $\QA$ and $\QB$ on  \cref{ln:preproc-begin-1,ln:preproc-mid-1,ln:preproc-end-1,ln:recovery} can be done without explicitly forming $\QA$ and $\QB$ by using \verb|xORMQR|, in the real case, or \verb|xUNMQR|, in the complex case.
  These two functions also require only $\max(m,n)$ flovals of work memory in precision $\uhigh$.
  Finally, the solution of the triangular Sylvester equation on \cref{ln:initial-sylv-1} can be performed with \verb|xTRSYL|.
  As this routing overwrites the right-hand side of the equation with the solution, a copy of $F$ must be computed, for an additional $mn$ flovals in precision $\uhigh$.
  \funstit{} requires another $mn$ flovals in precision $\uhigh$, thus overall \cref{alg:sylv-mp} requires least $2\max(m,n)^{2} + 3\max(m,n)$ flovals in precision $\ulow$ and \mbox{$2(m^{2} + n^{2}) + 2mn + \max(m,n)$} flovals in precision $\uhigh$.
  We note that the flovals in precision $\ulow$ can be reused as $\uhigh$ working memory after the low-precision Schur decomposition have been converted to high precision.

  For comparison, the standard Bartels--Stewart algorithm (see \cref{sec:BS-algorithm}) run in precision $\uhigh$ requires $m^2+n^2 + 3\max(m,n)$ flovals for the two Schur decompositions, where the Schur factors overwrite the input and an additional storage of $m^2+n^2$ is required for storing the unitary factors.
  The subsequent triangular Sylvester equation requires $mn$ flovals in precision $\uhigh$ for storing the right-hand side coefficient matrix, which is then overwritten by the solution $Y$. The final transformation of $Y$ to $X$ requires no extra memory.
  Therefore, the standard Bartels--Stewart algorithm requires at least $m^{2} + n^{2} + mn + 3\max(m,n)$ flovals in precision $\uhigh$. This means that the additional storage space required by \cref{alg:sylv-mp} is $2\max(m,n)^{2} + 3\max(m,n)$ flovals in precision $\ulow$ and $m^{2} + n^{2} + mn - 2\max(m,n)$ flovals in precision $\uhigh$.

  \subsection{Inversion of the unitary factors in high precision}
  \label{sec:mixed-precision-inv}

  \newcommand{\MA}{\ensuremath{M_{A}}}
  \newcommand{\MB}{\ensuremath{M_{B}}}
  As discussed in \cref{sec:mixed-precision-orth}, one cannot simply recover the solution $X$ to~\eqref{eq:sylv} in precision
  $\uhigh$ by inverting the two nearly-unitary matrices $\wUA^*$ and
  $\wUB$, which are unitary in precision $\ulow$.

  Recall from~\eqref{eq:rw-pert} that $\wDA = A - \wUA \wTA \wUA^*$ and $\wDB = B - \wUB \wTB \wUB^*$.
  If we write the matrix equation~\eqref{eq:sylv} as
  \begin{equation*}
    \wUA\bigl(\wTA + \wUA^{-1} \wDA \wUA^{-*}\bigr) \wUA^{*} X +
    X \wUB \bigl(\wTB + \wUB^{-1} \wDB \wUB^{-*}\bigr)\wUB^{*} =
    C,
  \end{equation*}
  and we multiply it by $\wUA^{-1}$ on the left and by $\wUB^{-*}$ on
  the right, we obtain
  \begin{equation}
    \label{eq:mult-inv}
    \bigl(\wTA + \wUA^{-1} \wDA \wUA^{-*}\bigr) \wUA^{*} X \wUB^{-*} +
    \wUA^{-1} X \wUB \bigl(\wTB + \wUB^{-1} \wDB \wUB^{-*}\bigr) =
    \wUA^{-1} C \wUB^{-*}.
  \end{equation}
  It is clear that we cannot substitute the two unknowns with a same
  matrix $Y$, because $\wUA^{*} X \wUB^{-*} \approx \wUA^{-1} X \wUB$
  only to precision $\ulow$.
  If, however, we set $Y:=\wUA^{*} X\wUB$, we can
  rewrite~\eqref{eq:mult-inv} as
  \begin{equation*}
    \bigl(\wTA + \wUA^{-1} \wDA \wUA^{-*}\bigr) Y\wUB^{-1} \wUB^{-*} +
    \wUA^{-1} \wUA^{-*} Y \bigl(\wTB + \wUB^{-1} \wDB \wUB^{-*}\bigr) =
    \wUA^{-1} C \wUB^{-*}.
  \end{equation*}
  Multiplying this expression by $\wUA^{*}\wUA$ on the left and by
  $\wUB^{*}\wUB$ on the right, we obtain
  \begin{equation}\label{eq:1}
    \wUA^{*}\wUA\bigl(\wTA + \wUA^{-1} \wDA \wUA^{-*}\bigr) Y +
    Y \bigl(\wTB + \wUB^{-1} \wDB \wUB^{-*}\bigr) \wUB^{*}\wUB =
    \wUA^{*} C \wUB.
  \end{equation}
  If we expand further, rewrite $\wUA^{*}\wUA= I +(\wUA^{*}\wUA-I)$ and
  $\wUB^{*}\wUB = I + (\wUB^{*}\wUB -I)$, and simplify like terms, we get
  \begin{equation}
    \label{eq:equation-with-ma-mb}
    (\wTA + \DTA)  Y + Y (\wTB + \DTB)  = \wUA^{*} C \wUB,
  \end{equation}
  where the two matrices
  \begin{equation}
    \label{eq:ma-mb}
    \begin{aligned}
      \DTA &= \wUA^{-1} \wDA \wUA^{-*} +
      (\wUA^{*}\wUA-I)
      \bigl(\wTA + \wUA^{-1} \wDA \wUA^{-*}\bigr),\\
      \DTB &= \wUB^{-1} \wDB \wUB^{-*} + \bigl(\wTB + \wUB^{-1} \wDB
      \wUB^{-*}\bigr)
      (\wUB^{*}\wUB-I)
    \end{aligned}
  \end{equation}
  are small entry-wise.
  By expanding the parentheses in~\eqref{eq:ma-mb}, we obtain
  \begin{align*}
    \DTA
    &= \wUA^{-1} \wDA \wUA^{-*} +
    \wUA^{*}\wUA\wTA + \wUA^{*} \wDA \wUA^{-*}
    - \wTA - \wUA^{-1} \wDA \wUA^{-*}\\
    &= \wUA^{*}\wUA\wTA + \wUA^{*} \wDA \wUA^{-*} - \wTA\\
    &= \wUA^{*}\wUA\wTA + \wUA^{*} (A - \wUA\wTA\wUA^{*}) \wUA^{-*} -
    \wTA\\
    &= \wUA^{*} A \wUA^{-*} - \wTA,
  \end{align*}
  and similarly $ \DTB = \wUB^{-1} B \wUB - \wTB.$

  \begin{algorithm2e}[t]
    \caption{Inversion-based mixed-precision Sylvester solver.}
    \label{alg:sylv-mp-inversion}
    \KwIn{$A\in\Cmm$, $B\in\Cnn$, $C\in\Cmn$, $\ulow$ and $\uhigh$ as in~\eqref{eq:prec-inequality}.}
    \KwOut{$X \in \Cmn$ such that $A X + X B \approx C$.}
    Compute Schur decomposition $A =: \wUA \wTA \wUA^*$ in precision
    $\ulow$.\;\label{ln:schurA-mp-inv}
    Compute Schur decomposition $B =: \wUB \wTB \wUB^*$ in precision
    $\ulow$.\;\label{ln:schurB-mp-inv}
    Compute LU decomposition of $\wUA^{*}$ in precision
    $\uhigh$.\;\label{ln:luA-mp-inv}
    Compute LU decomposition of $\wUB$ in precision
    $\uhigh$.\;\label{ln:luB-mp-inv}
    $F \gets \wUA^{*}C\wUB $, computed in precision
    $\uhigh$.\;\label{ln:preproc-begin-inv-1}
    $\LA \gets \wUA^{*} A \wUA^{-*} - \wTA$, computed in precision
    $\uhigh$.\;\label{ln:preproc-mid-inv-1}
    $\LB \gets \wUB^{-1} B \wUB - \wTB$, computed in precision
    $\uhigh$.\;\label{ln:preproc-end-inv-1}
    Find $Y_0$ such that $\wTA Y_0 + Y_0 \wTB = F$ in precision
    $\ulow$.\;\label{ln:first-sylv-solv-inv}
    $Y \gets$ \funstit($\wTA$, $\LA$, $\wTB$, $\LB$, $F$,
    $Y_0$) in precision $\uhigh$.\;\label{ln:call-inv-stit}
    $X \gets \wUA^{-*} Y \wUB^{-1}$ in precision
    $\uhigh$.\;\label{ln:recovery-inv}
  \end{algorithm2e}

  If we use \cref{alg:pert-sylv-tri-stationary} on~\eqref{eq:equation-with-ma-mb},
  we obtain \cref{alg:sylv-mp-inversion}.
  For the sufficient condition~\eqref{eq:norm-2-fro},
  it can be shown that
  \begin{align*}
    \norm{\DTA}
    &= \norm{\wUA^{*} A \wUA^{-*} - \wTA}
    = \norm{\wUA^{*}(\wUA \wTA \wUA^{*} + \wDA) \wUA^{-*} - \wTA}\\
    &= \norm{(\wUA^{*}\wUA - I) \wTA + \wUA^{*}\wDA\wUA^{-*}}
    \le \norm{\wUA^{*}\wUA - I}\norm{\wTA} +
    \norm{\wUA^{*}}\norm{\wUA^{-*}}\norm{\wDA}\\
    &\approx \ulow\bigl(\norm{\wTA} + \kappa(\wUA^{*})\norm{A}\bigr),
  \end{align*}
  and similarly that $\norm{\DTB} \lesssim \ulow\bigl(\norm{\wTB} +
  \kappa(\wUB)\norm{B}\bigr)$.

  \subsubsection{Computational cost and storage requirement}
  \label{sec:mixed-precision-inv-cost}
  We now discuss the computational cost of this algorithm, which is also reported
  in \cref{tab:mp-comp-cost} for ease of comparison.
  Computing the Schur decompositions of $A$ and $B$ on
  \cref{ln:schurA-mp-inv,ln:schurB-mp-inv} requires
  $25(m^3+n^3)$ flops in precision $\ulow$, whereas computing the LU
  decomposition of $\UA^{*}$ and $\UB$ on
  \cref{ln:luA-mp-inv,ln:luB-mp-inv}, which will be used to
  solve the linear systems later on, requires $\frac{2}{3}(m^3+n^3)$
  flops in precision $\uhigh$.
  The preprocessing step on \cref{ln:recovery-inv} requires two
  matrix products, which overall account for an additional $2mn(m+n)$
  flops in precision $\uhigh$.

  Computing $\LA$, $\LB$, and $F$ on \cref{ln:preproc-begin-inv-1,ln:preproc-mid-inv-1,ln:preproc-end-inv-1} requires
  \mbox{$4(m^3+n^3) + 2mn(m+n)$} flops in precision $\uhigh$, and computing $Y_0$ requires $mn(m+n)$ flops in precision~$\ulow$.
  The call to \funstit{} on \cref{ln:call-inv-stit} requires
  $3imn(m+n)$ flops in precision $\uhigh$, where $i$ is the number of
  iterations requires by the function that solves the triangular Sylvester
  equation.
  In this case, \cref{alg:sylv-mp-inversion} can solve~\eqref{eq:sylv}
  with $25(m^3+n^3) + mn(m+n)$ low-precision flops and
  $\big(4+\frac{2}{3}\big)(m^3+n^3)+(4+3i)mn(m+n)$ high-precision flops and~\eqref{eq:lyap} with $27n^{3}$ low-precision and
  $(12+\frac{2}{3}+6i)n^{3}$ high-precision flops.

  The analysis of the additional memory needed is similar to that in \cref{sec:mixed-precision-orth-cost}.
  Computing the Schur factorizations in low precision and converting them to high precision requires $\max(m,n)^{2} + 3 \max(m,n)$ floval in precision $\ulow$ and $2 (m^{2} + n^{2})$ floval in precision $\uhigh$.

  To minimize the additional memory required to compute $F$, $\LA$, and $\LB$, care in the evaluation order is needed, because \verb|xGEMM| cannot accumulate the result of a matrix product on either of its two factors.
  A possible solution comprises two steps.
  First, we allocate $\max(m,n)^{2}$ flovals to a matrix $Z$ and then compute

\algorithmstyle{plain}
\begin{algorithm2e}[H]
$Z \gets \wUA^* C$\;
$C \gets Z \wUB$  \tcc*{$C$ now contains the updated right-hand side}
$Z \gets \wUA^{*} A$\;
$A \gets Z$\;
$Z \gets B \wUB$\;
$B \gets Z$\;
\end{algorithm2e}
\algorithmstyle{ruled}

At this point, the matrices $\wUA$ and $\wUB$ cease to be needed, and we can use \verb|xGETRF| to compute their LU decompositions.
This function modifies the input but does not require any work memory.
Once the LU decompositions have been computed, we can use \verb|xTRSM|, a BLAS routing that does not require additional work memory, to update $\LA$ and $\LB$ without additional memory required.
As already discussed in the case of \cref{alg:sylv-mp-inversion}, the solution of the equation on \cref{ln:first-sylv-solv-inv} does not require any additional memory, while the iteration on \cref{ln:call-inv-stit} requires $mn$ flovals in precision $\uhigh$, for which the matrix $Z$ can be used.

Therefore, \cref{alg:sylv-mp-inversion} requires at least $2\max(m,n)^{2} + 3 \max(m,n)$ in precision $\ulow$ and $2(m^{2} + n^{2} + mn) + \max(m,n)$ flovals in precision $\uhigh$.
Compared with the Bartels--Steward algoirthm, \cref{alg:sylv-mp-inversion} requires additionally $2\max(m,n)^{2} + 3\max(m,n)$ flovals in precision $\ulow$ and $m^{2} + n^{2} + mn - 2\max(m,n)$ flovals in precision $\uhigh$.

  \newcommand{\ratio}{\ensuremath{\rho}}
  \newcommand{\funk}[3][(\rho)]{\ensuremath{\varphi_{#2}^{#3}#1}}
  \newcommand{\cratio}[3][(\rho)]{\ensuremath{C_{#2}^{#3}#1}}

  \section{A flop-based computational model}
  \label{sec:flop-based-comp-model}
  When will the computational cost of the mixed-precision algorithms be lower than that of the Bartels--Stewart algorithm in high precision?
  We answer this question with the following computational model.

  Let $\ratio$ be the ratio of the computational cost of a flop in
  precision~$\ulow$ to one in precision~$\uhigh$.
  Computing in precision $\uhigh$ is usually more expensive than computing in
  precision~$\ulow$, so we should expect $\ratio \le 1$ and in practice
  $\ratio \ll 1$ in most cases of practical interest.
  For each $\ratio$, we can find the maximum number of iterations that \cref{alg:sylv-mp,alg:sylv-mp-inversion} can perform while asymptotically requiring fewer operations than the Bartels--Stewart approach in high precision.

  To do this, we define, for each algorithm, a function of $\ratio$ that represents the iteration threshold.
  For \cref{alg:sylv-mp}, we have
  \begin{equation}
    \label{eq:funk-orth}
    \funk{S}{1} = \frac{(19-25\ratio)(m^{3}+n^{3}) +
      (1-\ratio)mn(m+n)}{3mn(m+n)},
  \end{equation}
  for the Sylvester equation, and
  \begin{equation}
    \label{eq:funk-orth-lyap}
    \funk{L}{1} =
    \frac{21-27\ratio}{6},
  \end{equation}
  for the Lyapunov equation.
  For \cref{alg:sylv-mp-inversion}, we obtain
  \begin{equation}
    \label{eq:funk-inv}
    \funk{S}{2} =
    \dfrac{\big(20 + \frac{1}{3}-25\ratio\big)(m^{3}+n^{3}) +
      (1-\ratio)mn(m+n)}{3mn(m+n)},
  \end{equation}
  for the Sylvester equation, and
  \begin{equation}
    \label{eq:funk-inv-lyap}
    \funk{L}{2} =
    \frac{22+\frac{1}{3}-27\ratio{}}{6},
  \end{equation}
  for the Lyapunov equation.

\input{figs/optk-funk.tex}

  In the left panel of \cref{fig:optk}, we plot these four quantities, for $m=n$, and the corresponding integer parts, for values of $\rho$ between 0 (flops in precision \ulow{} have no cost) and 1 (flops in precision \ulow{} have the same cost as those in precision \uhigh{}).

  The results suggest that, for Sylvester equations, the mixed-precision approach can be computationally advantageous for $\rho$ as large as 0.7, as long as one iteration is sufficient to achieve convergence to the desired accuracy; for lower values of $\rho$, convergence in up to 7 iterations can bring potential performance benefits.
  For Lyapunov equations, less work can be performed in low precision, and even though potential gains are possible for $\rho$ as large as 0.6, to keep the computational cost of the mixed-precision algorithm below that of the purely high-precision alternative, we cannot afford more than 3 iterations.

  We also note that the curves for \cref{alg:sylv-mp-inversion} are slightly more favorable than those for \cref{alg:sylv-mp}.
  This is because the former has lower computational cost, requiring $\big(1+\frac{1}{3}\big)(m^3+n^3)$ fewer high-precision flops during the pre-processing stage.

  With our computational model, we can gauge the cost ratios of the mixed-precision algorithms relative to the high-precision Bartels--Stewart algorithm for a given number $i$ of refinement steps.
  In the right panel of \cref{fig:optk}, we consider the case $i = 1$, which is the ideal scenario for the mixed-precision algorithms.
  Using the flop counts in \cref{sec:mixed-precision-orth-cost,sec:mixed-precision-inv-cost}, we compute these ratios, which we denote by $\cratio{X}{Y}$, where $X$ is $S$ for Sylvester or $L$ for Laypunov, and $Y$ is $1$ for \cref{alg:sylv-mp} and $2$ for \cref{alg:sylv-mp-inversion}.

  The results suggest that, for Sylvester equations, the mixed-precision algorithms can reduce the computational cost by up to 60\%, when the cost of low-precision flops becomes negligible. When $\rho=0.5$, we can expect the mixed-precision algorithms to be 20\% cheaper than the Bartels--Stewart algorithm. For Lyapunov equations, on the other hand, the computational cost can be reduced by 40\% at most, and the savings become marginal for $\rho=0.5$.

  \section{GMRES-IR for the Sylvester matrix equation}\label{sec:gmres-ir}

  It is possible to solve the Sylvester equation~\eqref{eq:sylv} by applying an iterative algorithm to the equivalent formulation~\eqref{eq:sylv-ls}.
  In a mixed-precision setting, one might derive a variant of the GMRES-IR algorithms~\cite{abhl24},~\cite{cahi17},~\cite{cahi18} tailored to the Sylvester equation.
  Good preconditioners are often necessary for such algorithms to be efficient.

  In the setting where the Schur decompositions $A = \UA \TA \UA^*$ and $B = \UB \TB \UB^*$ are computed in low precision, one can implicitly apply the obvious preconditioner $M_f^{-1}=(I_n\otimes A+B^T\otimes I_m)^{-1}$, which can yield an efficient algorithm if the GMRES solver converge sufficiently quickly.
  However, since applying the preconditioner involves solving Sylvester equations---which must necessarily be done in a precision lower than the working precision for this approach to be computationally sensible---the overall Schur-preconditioned GMRES-IR algorithm only converges for problems whose condition numbers are well bounded, depending on the unit roundoff of the precision at which the preconditioner is formed and applied.
  A detailed discussion is beyond the scope of this manuscript.

  \section{Numerical experiments}\label{sec:numeric.experiments}

  We compare \cref{alg:sylv-mp,alg:sylv-mp-inversion} with
  the algorithm of Bartels and Stewart~\cite{bast72} run entirely in high precision.
  The experiments were run using the GNU/Linux release of MATLAB 9.14.0 (R2023a Update 7) on a machine equipped with a 32-core AMD EPYC 9354P CPU.
  The high precision was set to binary64, while precisions lower than binary64 were simulated using the CPFloat library~\cite{fami23}.
  The code to repeat our experiments is available on GitHub.\footnote{\url{https://github.com/north-numerical-computing/mixed-precision-sylvester}}

  \input{figs/mixedprecision.tex}

  We test the mixed-precision algorithms on matrix equations arising from various applications. Our test set contains 19 Sylvester and 12 Lyapunov equations from the literature~\cite{bai11},~\cite{blw07},~\cite{benn04},~\cite{bmsv99},~\cite{bequ05},~\cite{bqq05},~\cite{besa05},~\cite{dfg22},~\cite{hbz05},~\cite{hure92},~\cite{lzz20},~\cite{wlm13},~\cite{zwt15}; the coefficients of these equations have order between 31 and 1,668, with the majority being in the few hundreds.
  The accuracy is gauged by evaluating in binary64 arithmetic the residual~\eqref{eq:sylv-norm-res} in the Frobenius norm.

  We compare the performance of the following codes:
  \begin{itemize}
    \item \lyap, the built-in MATLAB function \verb#lyap#, which calls the built-in function \sylvester{} if $B \neq A^*$;
    \item \mporth, a MATLAB implementation of
    \cref{alg:sylv-mp}; and
    \item \mpinv, a MATLAB implementation of
    \cref{alg:sylv-mp-inversion}.
  \end{itemize}

 In \cref{alg:pert-sylv-tri-stationary}, we allow a maximum of 20 iterations and check for convergence in the Frobenius norm, setting \mbox{$\varepsilon=10^{-12} \max(m,n)$}.

\subsection{TensorFloat-32 as low precision}
\label{sec:tensorfloat-32}

  In the first experiment, we use simulated TensorFloat-32 as low precision. The results are presented in \cref{fig:pert-mp}, where the matrices are sorted from by decreasing values of $\kappa_{\infty}(M_{f})$, for $M_{f}$ in~\eqref{eq:sylv-ls}.

  Overall, \cref{alg:sylv-mp,alg:sylv-mp-inversion} deliver accuracy comparable to that of \lyap{} and \sylvester{}.
  The built-in function \sylvester{} exhibits some instability for Sylvester equations 11 and 12.
  Both mixed-precision algorithms fail on Sylvester equation 1. For this problem, $\lambda\in\Lambda(\wTA)$ and $-\lambda\in\Lambda(\wTB)$ for some $\lambda\in\mathbb{R}$, thus $0\in\Lambda(M_{f})$ and the corresponding Sylvester equation is singular.
  Thus, the triangular solves on \cref{ln:initial-sylv-1} of \cref{alg:sylv-mp} and \cref{ln:first-sylv-solv-inv} of \cref{alg:sylv-mp-inversion} produce a matrix containing NaNs.
  The convergence of the mixed-precision algorithms is also slow for the Lyapunov equations~3 and~4.
  This behavior is consistent with our analysis: as reported in \cref{tab:ir-linear-system}, using TensorFloat-32 as low precision only guarantees convergence for equations whose condition number is of order $10^{4}$ or less, which is only true for Sylvester equations~7--19 and Lyapunov equations 10--12.

  We also note that there is a correlation between the condition number of the matrix equation and the number of refinement steps needed.
  Following the cost analysis in \cref{sec:mixed-precision}, we can conclude that \mporth{} and \mpinv{} will require fewer flops than the Bartels--Stewart algorithm in most test cases if $\rho$ for the pair \mbox{TensorFloat-32/binary64} is at most 0.16 and 0.32, for five and four refinement steps, respectively.

  To assess how realistic this condition is, consider the throughput of existing hardware that supports mixed-precision matrix operations.
  On the latest NVIDIA Grace--Blackwell GB200 and GB300 superchips, if tensor cores are enabled, the peak throughput of dense matrix--matrix multiplication in TensorFloat-32 is 90 petaflops/s.
  Using binary64 tensor cores, the same kernel only achieves a throughput of 0.1 petaflop/s on the GB200 and 2.88 petaflop/s on the GB300~\cite{nv25}.
  Therefore, the speedup of TensorFloat-32 over binary64 is about $900\times$ for the GB200 and $31\times$ for the GB300.

  For $m=n$ and $\rho\le 0.1$, \cref{alg:sylv-mp} with four refinement steps and \cref{alg:sylv-mp-inversion} with five steps require the equivalent of $\big(55+\frac{1}{5}\big)n^3$ and $\big(46+\frac{8}{15}\big) n^3$ binary64 flops, respectively, for the Sylvester equations.
  This computational cost is at least 9.2\% and 22.4\% lower than that of the binary64 Bartels--Stewart algorithm, which requires $60n^3$ flops.

  The Lyapunov equations in our test set are rather ill conditioned, and three refinement steps are generally insufficient for the mixed-precision algorithms to converge.
  Therefore, we should expect our mixed-precision algorithms to be slower than the Bartels--Stewart algorithm.

  \subsection{Custom low-precision format with 16-bit significand}

  \input{figs/mixedprecision2}

  Now we consider a custom 24-bit (3-byte) floating-point format that has the same exponent range as TensorFloat-32 but 16 rather than 11 significant bits.
  This format has roughly the same dynamic range as TensorFloat-32 but is more accurate by a factor $2^5=32$.

  We repeat the experiment in~\cref{sec:tensorfloat-32} using this custom format as low precision. The results are reported in \cref{fig:pert-mp2}.
  Decreasing the unit roundoff of the low precision has cured the instability of \mporth{} and \mpinv{}.
  Both mixed-precision algorithms converge in 2 to 3 iterations in most cases, which suggests that they will be faster than the Bartels--Stewart algorithm if $\rho\lesssim 0.4$ for the Sylvester equations and $\rho\lesssim0.1$ for the Lyapunov equations.

  This experiment also shows that increasing $\ulow$, and therefore reducing the precision, does not necessarily improve the time-to-solution of the mixed-precision algorithm.
  In fact, lower precision will speed up the computation of the Schur decompositions, but it may increase the number of iterations required, leading to a longer runtime overall.

  \section{Conclusions}\label{sec:conclusions}
  We have derived two algorithms for solving the Sylvester equation using two floating-point precisions.
  The main building block is the stationary iteration in \cref{alg:pert-sylv-tri-stationary}, which iteratively refines a solution to the perturbed quasi-triangular Sylvester equation.
  We have analyzed the convergence of this method and its attainable residual in a two-precision setting.

  This iteration can be used, for example, to refine an approximate solution obtained using the quasi-triangular factors of low-precision Schur decompositions, so that it is accurate to high precision.
  This observation has allowed us to develop two new approaches to solve the Sylvester and Lyapunov equation in two precisions.
  The two algorithms we have proposed leverage either orthonormalization (\cref{alg:sylv-mp}) or inversion (\cref{alg:sylv-mp-inversion}) to obtain high-precision unitary matrices that can be used to recover the accurate full solution from the accurate quasi-triangular solution obtained via iterative refinement.

  We have proposed a model to compare the flop count of mixed-precision and mono-precision algorithms.
  This model is then employed to compare the computational cost of our new algorithms with the cost of the Bartels--Stewart algorithm run entirely in high precision.
  Our numerical experiments, run on Sylvester and Lyapunov equations from the literature, show that the accuracy of the new approaches is comparable with that of the high-precision Bartels--Stewart algorithm.
  The experiments also suggest that, for Sylvester equations, a performance gain should be expected from high-performance implementations on existing hardware.
  Such implementations should target hardware for which a low-precision implementation of the QR algorithm is available, so that the performance gain of low precision can be assessed on real hardware.
  For Lyapunov equations, on the other hand, our experiments suggests that a performance gain should only be expected for well-conditioned equations, where one or two iterations are sufficient to achieve convergence.

  An open question is the performance of the mixed-precision algorithms on Sylvester equations with unbalanced coefficients ($m\gg n$ or $n\gg m$), where the high-precision stationary iteration in \cref{alg:pert-sylv-tri-stationary} becomes relatively much cheaper, potentially changing the overall cost landscape.
  We will also consider whether extending our approach to more than two precisions has the potential for further acceleration.
  Finally, we will examine whether the results discussed here can be applied to T- and $\star$-Sylvester matrix equations~\cite{dgks16}, as well as to systems consisting of T-, $\star$-, and Sylvester matrix equations~\cite{dmka15},~\cite{joka02}.

\section*{Acknowledgments} The experimental work was undertaken on the Aire HPC system at the University of Leeds, UK. The authors thank the four anonymous reviewers for their comments on an earlier draft of this manuscript.

  \bibliographystyle{siamplain-doi}
  \bibliography{references}

\end{document}

%% file: tikz.tex
\newlength\figwidth
\newlength\figheight
\pgfplotsset{
  legend style={font=\footnotesize},
  execute at begin axis={
    \pgfplotsset{
      axis line style = thick,
      tick label style={font=\normalsize},
      ylabel style={rotate=-90},
      xticklabel style={
        /pgf/number format/fixed,
        /pgf/number format/precision=5
      },
      scale only axis,
      scaled x ticks=false,
      yminorticks=true,
      axis background/.style = {fill=white},
      legend cell align = {left},
      legend style={thick,mark size=3pt},
    }
  },
  grid=both,
  grid style={line width=.75pt, draw=gray!50,densely dotted},
  major grid style={line width=.9pt,draw=gray!90},
  max space between ticks=30pt
}

%% file: tikz_lists.tex
\pgfplotscreateplotcyclelist{sylvlist}{
  every mark/.append style={solid}, line width=1.5pt, mark size=3pt, mark=o, only marks, Blue\\
  every mark/.append style={solid}, line width=1.5pt, mark size=3pt, mark=star, only marks, BrickRed\\
  every mark/.append style={solid}, line width=1.5pt, mark size=3pt, mark=diamond, only marks, JungleGreen\\
  every mark/.append style={solid}, line width=1.5pt, mark size=3pt, mark=none, solid, Cyan\\
}

\pgfplotscreateplotcyclelist{mplist}{
    ultra thick,mark=none, BurntOrange\\
    every mark/.append style={solid}, line width=1.5pt, mark size=3pt, mark=o, only marks, BrickRed\\
    every mark/.append style={solid}, line width=1.5pt, mark size=3pt, mark=triangle, only marks, JungleGreen\\
    every mark/.append style={solid}, line width=1.5pt, mark size=3pt, mark=square, only marks, PineGreen\\
    every mark/.append style={solid}, line width=1.5pt, mark size=3pt, mark=star, only marks, Red\\
    every mark/.append style={solid}, line width=1.5pt, mark size=3pt, mark=asterisk, only marks, Orange\\
  }

%% file: metadata.tex
\newcommand{\thetitle}{Mixed-precision algorithms for solving\\
  the Sylvester matrix equation}
\newcommand{\thetitlenote}{Version of \today.}
\newcommand{\thefunding}{The work of Andrii Dmytryshyn was supported by the Swedish Research Council (VR) under grant 2021-05393.}

\newcommand{\theabstract}{
We consider the solution of the Sylvester equation $AX+XB=C$ in mixed precision.
We derive a new iterative refinement scheme to solve perturbed quasi-triangular Sylvester equations; our rounding error analysis provides sufficient conditions for convergence and a bound on the attainable relative residual.
We leverage this iterative scheme to solve the general Sylvester equation.
The new algorithms compute the Schur decomposition of the coefficient matrices $A$ and $B$ in lower than working precision, use the low-precision Schur factors to obtain an approximate solution to the perturbed quasi-triangular equation, and iteratively refine it to obtain a working-precision solution.
In order to solve the original equation to working precision, the unitary Schur factors of the coefficient matrices must be unitary to working precision, but this is not the case if the Schur decomposition is computed in low precision.
We propose two effective approaches to address this: one is based on re-orthonormalization in working precision, and the other on explicit inversion of the almost-unitary factors.
The two mixed-precision algorithms thus obtained are tested on various Sylvester and Lyapunov equations from the literature.
Our numerical experiments show that, for both types of equations, the new algorithms are at least as accurate as existing ones.
Our cost analysis, on the other hand, suggests that they would typically be faster than mono-precision alternatives if implemented on hardware that natively supports low precision.}

\newcommand{\thekeywords}{Sylvester equation,
    iterative refinement,
    mixed-precision,
    rounding error analysis,
    Schur decomposition,
    orthonormalization
}

\newcommand{\themsc}{65F45,  
    15A24, 
    65F10, 
    65G50 
}

\newcommand{\theauthori}{Andrii Dmytryshyn}

\newcommand{\theaffiliationi}{
    Department of Mathematical Sciences,
    Chalmers University of Technology and University of Gothenburg,
    41296 Gothenburg,
    Sweden}
\newcommand{\theemaili}{andrii@chalmers.se}

\newcommand{\theauthorii}{Massimiliano Fasi}

\newcommand{\theaffiliationii}{
    School of Computer Science,
    University of Leeds,
    Woodhouse Lane, Leeds LS2 9JT, UK}
\newcommand{\theemailii}{m.fasi@leeds.ac.uk}

\newcommand{\theauthoriii}{Nicholas J.~Higham}

\newcommand{\theaffiliationiii}{
    The author is deceased. Former address:
    Department of Mathematics,
    University of Manchester,
    Oxford Road, Manchester M13 9PL, UK}

\newcommand{\theauthoriv}{Xiaobo Liu}

\newcommand{\theaffiliationiv}{
     Max Planck Institute for Dynamics of Complex Technical Systems,
     Sandtorstraße, 39106 Magdeburg, Germany}
\newcommand{\theemailiv}{xliu@mpi-magdeburg.mpg.de}

%% file: commands.tex
\DeclareMathOperator{\sep}{sep}
\DeclareMathOperator{\vvec}{vec}

\DeclareMathOperator{\cond}{cond}

\DeclarePairedDelimiter{\norm}{\lVert}{\rVert}

\DeclarePairedDelimiter{\normt}{\lVert}{\rVert_2}
\DeclarePairedDelimiter{\abs}{\lvert}{\rvert}
\DeclarePairedDelimiter{\normo}{\lVert}{\rVert_\infty}

\renewcommand{\Delta}{\varDelta}

\renewcommand{\le}{\leqslant}

\numberwithin{equation}{section}
\numberwithin{table}{section}
\numberwithin{figure}{section}

\DontPrintSemicolon
\SetKwProg{Function}{function}{}{}
\let\oldnl\nl

\newcommand{\nonl}{\renewcommand{\nl}{\let\nl\oldnl}} 
\SetKwComment{Comment}{}{}
\makeatletter
\newcommand{\algorithmstyle}[1]{\renewcommand{\algocf@style}{#1}}
\makeatother
\SetKwInOut{Parameter}{parameter}
\SetKwRepeat{Repeat}{repeat}{until}

\renewcommand{\Sigma}{\varSigma}

\definecolor{mygreen}{RGB}{28,172,0}
\definecolor{mylilas}{RGB}{170,55,241}

\newcommand{\ulow}{\ensuremath{u_\ell}}
\newcommand{\uhigh}{\ensuremath{u_h}}

\newcommand{\DA}{\ensuremath{\Delta A}}
\newcommand{\wDA}{\ensuremath{\wh{\DA}}}
\newcommand{\DB}{\ensuremath{\Delta B}}
\newcommand{\wDB}{\ensuremath{\wh{\DB}}}

\newcommand{\DM}{\ensuremath{\Delta M}}

\newcommand{\wh}{\widehat}
\newcommand{\wt}{\widetilde}
\newcommand{\Ctri}{\ensuremath{\wt C}}

\newcommand{\DT}{\ensuremath{\Delta T}}

\newcommand{\TA}{\ensuremath{T_{A}}}
\newcommand{\wTA}{\ensuremath{\wh{T}_{A}}}
\newcommand{\DTA}{\ensuremath{\DT_{A}}}
\newcommand{\UA}{\ensuremath{U_{A}}}
\newcommand{\wUA}{\ensuremath{\wh{U}_{A}}}
\newcommand{\EA}{\ensuremath{E_{A}}}
\newcommand{\QA}{\ensuremath{Q_{A}}}
\newcommand{\RA}{\ensuremath{R_{A}}}
\newcommand{\LA}{\ensuremath{L_{A}}}

\newcommand{\TB}{\ensuremath{T_{B}}}
\newcommand{\wTB}{\ensuremath{\wh{T}_{B}}}
\newcommand{\DTB}{\ensuremath{\DT_{B}}}
\newcommand{\UB}{\ensuremath{U_{B}}}
\newcommand{\wUB}{\ensuremath{\wh{U}_{B}}}
\newcommand{\EB}{\ensuremath{E_{B}}}
\newcommand{\QB}{\ensuremath{Q_{B}}}
\newcommand{\RB}{\ensuremath{R_{B}}}
\newcommand{\LB}{\ensuremath{L_{B}}}

\newcommand{\mathset}[1]{\ensuremath{\mathbb{#1}}}

\newcommand{\kbyk}{\ensuremath{k \times k}}
\newcommand{\nbyn}{\ensuremath{n \times n}}
\newcommand{\mbym}{\ensuremath{m \times m}}
\newcommand{\mbyn}{\ensuremath{m \times n}}

\newcommand{\Ckk}{\mathset{C}^{\kbyk}}
\newcommand{\Cnn}{\mathset{C}^{\nbyn}}

\newcommand{\Cn}{\mathset{C}^{n}}
\newcommand{\Cmm}{\mathset{C}^{\mbym}}
\newcommand{\Cmn}{\mathset{C}^{\mbyn}}
\newcommand{\Css}{\mathset{C}^{s\times s}}
\newcommand{\Cs}{\mathset{C}^{s}}

\newcommand{\myalpha}{\alpha}
\newcommand{\mybeta}{\beta}
\newcommand{\mygamma}{\gamma}

\newcommand{\lyap}{\texttt{lyap}}
\newcommand{\sylvester}{\texttt{sylvester}}
\newcommand{\mporth}{\texttt{mp\_orth}}
\newcommand{\mpinv}{\texttt{mp\_inv}}

\newcommand{\whd}{\ensuremath{\wh{d}}}
\newcommand{\whr}{\ensuremath{\wh{r}}}
\usepackage{tcolorbox}
\tcbuselibrary{skins}
\tcbuselibrary{breakable}
\tcbset{shield externalize}

%% file: figs/optk-funk.tex
\pgfplotscreateplotcyclelist{method-list}{
  mark=none, draw=BrickRed\\
  only marks, mark=*, mark options={fill=BrickRed!40, draw opacity=0}\\
  mark=none, densely dashed, draw=JungleGreen\\
  only marks, mark=asterisk, draw=JungleGreen\\
  mark=none, dashed, draw=RoyalBlue\\
  only marks, mark=square, draw=RoyalBlue\\
  mark=none, densely dotted, draw=Orange\\
  only marks, mark=diamond, draw=Orange\\
}

\pgfplotscreateplotcyclelist{method-list-2}{
  mark=none, densely dotted, draw=Orange\\
  mark=none, dashed, draw=RoyalBlue\\
  mark=none, densely dashed, draw=JungleGreen\\
  mark=none, draw=BrickRed\\
}

\newcommand{\drawoptk}[1]{\addplot+[
  mark repeat=4, mark phase=2,
  line width=1pt] table[x=rho, y=optk_stit_#1] {\dataoptk};}
\newcommand{\drawfunk}[1]{\addplot+[
  line width=1.8pt] table[x=rho, y=funk_stit_#1] {\dataoptk};}

\begin{figure}[t!]
  \pgfplotstableread{figs/optk.dat}\dataoptk
  \pgfplotstableread{figs/cost_ratio.dat}\datacost
  \crefname{algorithm}{Alg.}{Algs.}

  \centering

  \subfigure[Number of iterations.]{
  \begin{tikzpicture}[trim axis right, trim axis left]
    \begin{axis}[
      set layers,
      mark layer=axis grid,
      width=0.42\textwidth,
      height=0.34\textwidth,
      scale only axis,
      xmin=0, xmax=1,
      ymin=-3, ymax=9,
      xlabel={$\rho$},
      mark size=3pt,
      axis background/.style={fill=none},
      cycle list name = method-list,
      legend style={
        cells={mark size=3pt, line width=1.8pt},
        at={(0.9925,0.99)},
        anchor=north east,
        legend columns=2,
        legend cell align=left,
        /tikz/every even column/.append style={column sep=5pt},
      }
      ]
      \drawfunk{sylv_inv}\addlegendentry{$\funk[]{S}{2}$}
      \drawoptk{sylv_inv}\addlegendentry{{$\lfloor \funk{S}{2} \rfloor$}}
      \drawfunk{sylv_orth}\addlegendentry{$\funk[]{S}{1}$}
      \drawoptk{sylv_orth}\addlegendentry{$\lfloor \funk{S}{1} \rfloor$}

      \drawfunk{lyap_inv}\addlegendentry{$\funk[]{L}{2}$}
      \drawoptk{lyap_inv}\addlegendentry{$\lfloor \funk{L}{2} \rfloor$}
      \drawfunk{lyap_orth}\addlegendentry{$\funk[]{L}{1}$}
      \drawoptk{lyap_orth}\addlegendentry{$\lfloor \funk{L}{1} \rfloor$}
    \end{axis}
  \end{tikzpicture}}
  \hspace{20pt}
  \subfigure[Cost ratio.]{
  \begin{tikzpicture}[trim axis right, trim axis left]
    \begin{axis}[
      width = 0.42\textwidth,
      height = 0.34\textwidth,
      xmin = 0, xmax = 1,
      ymin = 0.3, ymax = 1.5,
      xlabel={$\rho$},
      legend pos=south west,
      mark size= 4pt,
      line width=1.8pt,
      cycle list name = method-list-2,
      legend style={
        cells={mark size=3pt, line width=1.8pt},
        at={(0.0075,0.99)},
        anchor=north west,
        legend columns=1,
        legend cell align=left,
        /tikz/every even column/.append style={column sep=5pt},
      }
      ]
      \pgfplotstableread{figs/cost_ratio.dat}\mydata
      \addplot table[x=rho, y=cost_stit_lyap_orth]{\mydata};\addlegendentry{\cratio{L}{1}}
      \addplot table[x=rho, y=cost_stit_lyap_inv]{\mydata};\addlegendentry{\cratio{L}{2}}
      \addplot table[x=rho, y=cost_stit_sylv_orth]{\mydata};\addlegendentry{\cratio{S}{1}}
      \addplot table[x=rho, y=cost_stit_sylv_inv]{\mydata};\addlegendentry{\cratio{S}{2}}

    \end{axis}
  \end{tikzpicture}}

  \caption{The panel on the left shows the maximum number of iterations for which the algorithms in \cref{sec:pert-sylv-algs} will be asymptotically faster than the Bartels--Stewart algorithm run in high precision.
  The plot shows the quantities in~\eqref{eq:funk-orth},~\eqref{eq:funk-orth-lyap},~\eqref{eq:funk-inv}, and~\eqref{eq:funk-inv-lyap} against the $\ratio$, the ratio of the computational cost of a low-precision to a high-precision operation.
  The panel on the right shows the cost ratio of mixed-precision iterative refinement for $i = 1$ to high-precision Bartels--Stewart.}
  \label{fig:optk}
\end{figure}

%% file: figs/mixedprecision.tex
\newcommand{\plotresidual}[3]{
  \begin{tikzpicture}[trim axis right, trim axis left]
    \begin{axis}[
      cycle list name = mplist,
      ymode=log, 
      ymin = #2, ymax = #3,
      xlabel={Matrix ID},
      ]
      \pgfplotstableread{#1}\mydata;
      \addplot table[x=id, y=condu]{\mydata};
      \addlegendentry{\;$\kappa_\infty(M_f) u$};
      \addplot table[x=id, y=res_sylv]{\mydata};
      \addlegendentry{\;\lyap}
      \addplot table[x=id, y=r_or]{\mydata};
      \addlegendentry{\;\mporth}
      \addplot table[x=id, y=r_in]{\mydata};
      \addlegendentry{\;\mpinv}
      \legend{};
    \end{axis}
  \end{tikzpicture}}

\newcommand{\plotresidualsylv}[3]{
  \begin{tikzpicture}[trim axis right, trim axis left]
    \begin{axis}[
      cycle list name = mplist,
      ymode=log, 
      ymin = #2, ymax = #3,
      xlabel={Matrix ID},
      ]
      \pgfplotstableread{#1}\mydata;
      \addplot table[x=id, y=condu]{\mydata};
      \addlegendentry{\;$\kappa_\infty(M_f) u$};
      \addplot table[x=id, y=res_sylv]{\mydata};
      \addlegendentry{\;\sylvester}
      \addplot table[x=id, y=r_or]{\mydata};
      \addlegendentry{\;\mporth}
      \addplot table[x=id, y=r_in]{\mydata};
      \addlegendentry{\;\mpinv}
      \legend{};
    \end{axis}
  \end{tikzpicture}}

\newcommand{\plotsteps}[3]{
  \begin{tikzpicture}[trim axis right, trim axis left]
    \begin{axis}[
      ybar, xtick align=inside, bar width = 0.3,
      ymin = #2, ymax = #3,
      xlabel={Matrix ID},
      ]
      \pgfplotstableread{#1}\mydata;
      \addplot[draw=black,fill=Green!50, pattern=north east lines, bar shift=-0.15]
      table[x=id, y=i_or]{\mydata};
      \addlegendentry{\;\mporth}
      \addplot[draw=black,fill=Blue!50, bar shift=0.15]
      table[x=id, y=i_in]{\mydata};
      \addlegendentry{\;\mpinv}
      \legend{};
    \end{axis}
\end{tikzpicture}}

\newcommand{\leftlegend}{
  \begin{tikzpicture}[trim axis right, trim axis left]
    \begin{axis}[
        height=0cm,
        width=0cm,
        cycle list name = mplist,
         /tikz/every even column/.append style={column sep=10pt},
         cells={line width=0.9pt, mark size=2.6},
        legend style={at={(0,0)}, anchor=center, legend columns=6, draw=none},
      ]
      \addplot(0,0);
      \addlegendentry{\;$\kappa_\infty(M_f) u$};
      \addplot(0,0);
      \addlegendentry{\;\lyap};
      \addplot(0,0);
      \addlegendentry{\;\mporth};
      \addlegendimage{draw=black,fill=Green!50, pattern=north east lines, ybar, ybar legend};
      \addlegendentry{\;\mporth};
      \addplot(0,0);
      \addlegendentry{\;\mpinv};
      \addlegendimage{draw=black,fill=Blue!50, bar shift=0.15, ybar, ybar legend};
      \addlegendentry{\;\mpinv};
    \end{axis}
  \end{tikzpicture}}

\begin{figure}[t!]
  \pgfplotsset{
    width = .4\linewidth,
    height = .3\linewidth,
  }
  \centering
  \pgfplotsset{xmin = 0, xmax = 20}
  \subfigure[Relative residual (Sylvester).]{
    \plotresidualsylv{figs/test_mixedprecision_sylv_tf32.dat}{1e-18}{1e1}}
  \hspace{20pt}
  \subfigure[Number of steps (Sylvester).]{
    \plotsteps{figs/test_mixedprecision_sylv_tf32.dat}{0}{20}}\\
  \pgfplotsset{xmin = 0, xmax = 13}
  \subfigure[Relative residual (Lyapunov).]{
    \plotresidual{figs/test_mixedprecision_lyap_tf32.dat}{1e-26}{1e1}}
  \hspace{20pt}
  \subfigure[Number of steps (Lyapunov).]{
    \plotsteps{figs/test_mixedprecision_lyap_tf32.dat}{0}{20}}\\
  \subfigure[Legend.]{
    \begin{minipage}{10cm}
      \centering
      \leftlegend{}
    \end{minipage}}
  \caption{Comparison of \texttt{lyap} and our MATLAB implementation of
    \cref{alg:sylv-mp,alg:sylv-mp-inversion} on matrix equations from the literature.
    The low-precision arithmetic is TensorFloat-32, for which $\ulow = 2^{-11}$.
    The top and bottom row refer to Sylvester and Lyapunov equations, respectively.
    Left: relative residual of the computed solution. 
    Right: number of iterative refinement steps.}
  \label{fig:pert-mp}
\end{figure}

%% file: figs/mixedprecision2.tex
\begin{figure}[t!]
  \pgfplotsset{
    width = .4\linewidth,
    height = .3\linewidth,
  }
  \centering
  \pgfplotsset{xmin = 0, xmax = 20}
  \subfigure[Relative residual (Sylvester).]{
    \plotresidualsylv{figs/test_mixedprecision_sylv_24bit.dat}{1e-18}{1e1}}
  \hspace{20pt}
  \subfigure[Number of steps (Sylvester).]{
    \plotsteps{figs/test_mixedprecision_sylv_24bit.dat}{0}{20}}\\
  \pgfplotsset{xmin = 0, xmax = 13}
  \subfigure[Relative residual (Lyapunov).]{
    \plotresidual{figs/test_mixedprecision_lyap_24bit.dat}{1e-26}{1e1}}
  \hspace{20pt}
  \subfigure[Number of steps (Lyapunov).]{
    \plotsteps{figs/test_mixedprecision_lyap_24bit.dat}{0}{20}}\\
  \subfigure[Legend.]{
    \begin{minipage}{10cm}
      \centering
      \leftlegend{}
  \end{minipage}}
  \caption{Comparison of \texttt{lyap} and our MATLAB implementation of
    \cref{alg:sylv-mp,alg:sylv-mp-inversion} on the same test set. The
    low-precision arithmetic is a custom 3-byte format, for which $\ulow =
    2^{-16}$.}
  \label{fig:pert-mp2}
\end{figure}

%% file: references.bib
@String{j-AMC                   = "Appl. Math. Comput."}

@String{j-AN                    = "Acta Numerica"}

@String{j-APP-MATH-COMP         = "Appl. Math. Comput."}

@String{j-APP-NUM-MATH          = "Appl. Numer. Math."}

@String{j-BIT                   = "BIT"}

@string{j-BLMS                  = "Bull. London Math. Soc."}

@String{j-CACM                  = "Comm. ACM"}

@String{j-COMP-MATH-APPL        = "Computers Math. Applic."}

@String{j-IEEE-AC               = "IEEE Trans. Automat. Control"}

@String{j-IEEE-TOAC             = "IEEE Trans. Autom. Control."}

@String{j-IMAJNA                = "IMA J. Numer. Anal."}

@String{j-IJHPCA                = "Int. J. High Perform. Comput. Appl."}

@String{j-JCM                   = "J. Comput. Math."}

@String{j-JSC                   = "J. Sci. Comput."}

@String{j-LAA                   = "Linear Algebra Appl."}

@String{j-MATH-COMP             = "Math. Comp."}

@String{j-NLAA                  = "Numer. Linear Algebra Appl."}

@String{j-SIMAX                 = "SIAM J. Matrix Anal. Appl."}

@String{j-SINUM                 = "SIAM J. Numer. Anal."}

@String{j-SIREV                 = "SIAM Rev."}

@String{j-SISC                  = "SIAM J. Sci. Comput."}

@String{j-SIAM                  = "J. Soc. Indust. Appl. Math"}

@String{j-TOMS                  = "ACM Trans. Math. Software"}

@String{pub-BH                  = "Birkh{\"{a}}user"}

@String{pub-BH:adr              = "Boston, MA, USA"}

@String{pub-JH                  = "Johns Hopkins University Press"}

@String{pub-JH:adr              = "Baltimore, MD, USA"}

@String{pub-SIAM                = "Society for Industrial and Applied
                                   Mathematics"}

@String{pub-SIAM:adr            = "Philadelphia, PA, USA"}

@String{pub-Springer            = "Spring{\-}er-Ver{\-}lag"}

@String{pub-Springer:adr-BH     = "Berlin, Heidelberg"}

@Article{dfg22,
  author = "Dmytryshyn, Andrii and Fasi, Massimiliano and Gulliksson, Mårten",
  title = "The Dynamical Functional Particle Method for Multi-Term Linear Matrix
                  Equations",
  journal = j-AMC,
  year = 2022,
  volume = 435,
  pages = 127458,
  month = dec,
  issn = "0096-3003",
  doi = "10.1016/j.amc.2022.127458"
  }

@Article{fami23,
  author = "Fasi, Massimiliano and Mikaitis, Mantas",
  title = "{CPFloat}: {A} {C} Library for Simulating Low-precision Arithmetic",
  journal = j-TOMS,
  year = 2023,
  volume = 49,
  number = 2,
  pages = "1-32",
  month = jun,
  issn = "1557-7295",
  doi = "10.1145/3585515"
  }

@Book{abbb99,
  author = "Anderson, E. and Bai, Z. and Bischof, C. and Blackford, L. S. and
                  Demmel, J. and Dongarra, J. and Du Croz, J. and Greenbaum, A.
                  and Hammarling, S. and McKenney, A. and Sorensen, D.",
  title = "{LAPACK} Users’ Guide",
  publisher = pub-SIAM,
  year = 1999,
  address = pub-SIAM:adr,
  month = jan,
  isbn = 9780898719604,
  doi = "10.1137/1.9780898719604"
  }

@Article{abhl24,
  author = "Amestoy, Patrick and Buttari, Alfredo and Higham, Nicholas J. and
                  L’Excellent, Jean-Yves and Mary, Theo and Vieublé, Bastien",
  title = "Five-Precision {GMRES}-Based Iterative Refinement",
  journal = j-SIMAX,
  year = 2024,
  volume = 45,
  number = 1,
  pages = "529–552",
  month = feb,
  issn = "1095-7162",
  doi = "10.1137/23m1549079"
  }

@Article{babe08,
  author = "Baur, Ulrike and Benner, Peter",
  title = "Gramian-Based Model Reduction for Data-Sparse Systems",
  journal = j-SISC,
  year = 2008,
  volume = 31,
  number = 1,
  pages = "776–798",
  month = jan,
  issn = "1095-7197",
  doi = "10.1137/070711578"
  }

@Article{bai11,
  author = "Bai, Zhongzhi",
  title = "On {Hermitian} and Skew-{Hermitian} Splitting Iteration Methods for
                  the Continuous {Sylvester} Equations",
  journal = j-JCM,
  year = 2011,
  volume = 29,
  number = 2,
  pages = "185-198",
  month = jun,
  issn = "1991-7139",
  doi = "10.4208/jcm.1009-m3152"
  }

@Article{bast72,
  author = "Richard H. Bartels and George W. Stewart",
  title = "Algorithm 432: {Solution} of the Matrix Equation {$AX+XB=C$}",
  journal = j-CACM,
  year = 1972,
  volume = 15,
  number = 9,
  pages = "820-826",
  doi = "10.1145/361573.361582"
  }

@Article{beme13,
  author = "Benner, Peter and Mena, Hermann",
  title = "Rosenbrock Methods for Solving {Riccati} Differential Equations",
  journal = j-IEEE-TOAC,
  year = 2013,
  volume = 58,
  number = 11,
  pages = "2950–2956",
  month = nov,
  issn = "1558-2523",
  doi = "10.1109/tac.2013.2258495"
  }

@InProceedings{benn04,
  author = "Benner, Peter",
  title = "Factorized Solution of {Sylvester} Equations with Applications in
                  Control",
  booktitle = "Proceedings of the 16th International Symposium on Mathematical
                  Theory of Networks and Systems",
  year = 2004,
  month = jul,
  address = "Leuven, Belgium",
  editors = "De Moor, B. and Motmans, B. and Willems, J. and Van Dooren, P. and
                  Blondel, V.",
  url = "https://mathweb.ucsd.edu/~helton/MTNSHISTORY/CONTENTS/2004LEUVEN/CDROM/papers/59.pdf"
  }

@InCollection{besa05,
  author = "Benner, Peter and Saak, Jens",
  title = "A Semi-Discretized Heat Transfer Model for Optimal Cooling of Steel
                  Profiles",
  booktitle = "Dimension Reduction of Large-Scale Systems",
  publisher = pub-Springer,
  year = 2005,
  editor = "Benner, Peter and Sorensen, Danny C. and Mehrmann, Volker",
  pages = "353-356",
  address = pub-Springer:adr-BH,
  journal = "Dimension Reduction of Large-Scale Systems",
  issn = "1439-7358",
  doi = "10.1007/3-540-27909-1_19",
  ISBN = 3540245456
  }

@Article{bhro97,
  author = "Bhatia, Rajendra and Rosenthal, Peter",
  title = "How and Why to Solve the Operator Equation {$AX-XB = Y$}",
  journal = j-BLMS,
  year = 1997,
  volume = 29,
  number = 1,
  pages = "1-21",
  month = jan,
  issn = "0024-6093",
  doi = "10.1112/s0024609396001828"
  }

@Article{bddd02,
  author = "Blackford, L. Susan and {et al.}",
  title = "An updated set of basic linear algebra subprograms ({BLAS})",
  journal = j-TOMS,
  year = 2002,
  volume = 28,
  number = 2,
  pages = "135–151",
  month = jun,
  issn = "1557-7295",
  doi = "10.1145/567806.567807"
  }

@Article{blw07,
  author = "Bao, Liang and Lin, Yiqin and Wei, Yimin",
  title = "A New Projection Method for Solving Large {Sylvester} Equations",
  journal = j-APP-NUM-MATH,
  year = 2007,
  volume = 57,
  number = "5-7",
  pages = "521-532",
  month = may,
  issn = "0168-9274",
  doi = "10.1016/j.apnum.2006.07.005",
  url = "https://doi.org/10.1016/j.apnum.2006.07.005",
  }

@InCollection{bmsv99,
  author = "Benner, Peter and Mehrmann, Volker and Sima, Vasile and Van Huffel,
                  Sabine and Varga, Andras",
  title = "SLICOT---A Subroutine Library in Systems and Control Theory",
  publisher = pub-BH,
  year = 1999,
  address = pub-BH:adr,
  pages = "499-539",
  booktitle = "Applied and Computational Control, Signals, and Circuits",
  editor = "Datta, Biswa N.",
  doi = "10.1007/978-1-4612-0571-5_10",
  ISBN = 9781461205715
  }

@Article{bqq05,
  author = "Benner, Peter and Quintana-Ortí, Enrique S. and Quintana-Ortí,
                  Gregorio",
  title = "Solving Stable {Sylvester} Equations via Rational Iterative Schemes",
  journal = j-JSC,
  year = 2005,
  volume = 28,
  number = 1,
  pages = "51-83",
  month = dec,
  issn = "1573-7691",
  doi = "10.1007/s10915-005-9007-2"
  }

@Article{care96,
  author = "Calvetti, D. and Reichel, L.",
  title = "Application of {ADI} Iterative Methods to the Restoration of Noisy
                  Images",
  journal = j-SIMAX,
  year = 1996,
  volume = 17,
  number = 1,
  pages = "165–186",
  month = jan,
  issn = "1095-7162",
  doi = "10.1137/s0895479894273687"
  }

@InCollection{bequ05,
  author = "Benner, Peter and Quintana-Ortí, Enrique S.",
  title = "Model Reduction Based on Spectral Projection Methods",
  booktitle = "Dimension Reduction of Large-Scale Systems",
  publisher = pub-Springer,
  year = 2005,
  editor = "Benner, Peter and Sorensen, Danny C. and Mehrmann, Volker",
  pages = "5-48",
  address = pub-Springer:adr-BH,
  issn = "1439-7358",
  doi = "10.1007/3-540-27909-1_1",
  ISBN = 9783540279099
  }

@Book{cofr03,
  author = "Martin J. Corless and Art Frazho",
  title = "Linear Systems and Control: An Operator Perspective",
  publisher = "CRC Press",
  year = 2003,
  address = "Boca Raton",
  pages = 390,
  doi = "10.1201/9780203911372",
  isbn = "978-0-203911-37-2"
  }

@Article{datt94,
  author = "Datta, Biswa Nath",
  title = "Linear and numerical linear algebra in control theory: some research
                  problems",
  journal = j-LAA,
  year = 1994,
  volume = "197–198",
  pages = "755–790",
  month = jan,
  issn = "0024-3795",
  doi = "10.1016/0024-3795(94)90512-6"
  }

@InCollection{hbz05,
  author = "Hohlfeld, Dennis and Bechtold, Tamara and Zappe, Hans",
  title = "Tunable Optical Filter",
  booktitle = "Dimension Reduction of Large-Scale Systems",
  publisher = pub-Springer,
  year = 2005,
  editor = "Benner, Peter and Sorensen, Danny C. and Mehrmann, Volker",
  pages = "337-340",
  address = pub-Springer:adr-BH,
  issn = "1439-7358",
  doi = "10.1007/3-540-27909-1_15",
  ISBN = 9783540279099
  }

@Article{cswz08,
  author = "Chen, Gang and Song, Yangqiu and Wang, Fei and Zhang, Changshui",
  title = "Semi-supervised Multi-label Learning by Solving a {Sylvester}
                  Equation",
  journal = "Proceedings of the 2008 SIAM International Conference on Data
                  Mining",
  year = 2008,
  month = apr,
  doi = "10.1137/1.9781611972788.37"
  }

@Article{dgks16,
  author = "Dopico, Froilán M. and González, Javier and Kressner, Daniel
                  and Simoncini, Valeria",
  title = "Projection methods for large-scale {$T$-Sylvester} equations",
  journal = j-MATH-COMP,
  year = 2016,
  volume = 85,
  number = 301,
  pages = "2427-2455",
  month = jan,
  issn = "1088-6842",
  doi = "10.1090/mcom/3081",
  url = "https://doi.org/10.1090/mcom/3081",
  }

@Article{dmka15,
  author = "Dmytryshyn, A. and Kågström, B.",
  title = "Coupled {Sylvester}-type Matrix Equations and Block Diagonalization",
  journal = j-SIMAX,
  year = 2015,
  volume = 36,
  number = 2,
  pages = "580-593",
  language = "en",
  issn = "0895-4798",
  doi = "10.1137/151005907"
  }

@Book{gova13,
  author = "Golub, Gene H. and Van Loan, Charles F.",
  title = "Matrix Computations",
  publisher = pub-JH,
  year = 2013,
  address = pub-JH:adr,
  edition = "4th",
  isbn = "1-4214-0794-9"
  }

@Article{hure92,
  author = "Hu, D.Y. and Reichel, L.",
  title = "{Krylov}-Subspace Methods for the {Sylvester} Equation",
  journal = j-LAA,
  year = 1992,
  volume = 172,
  pages = "283-313",
  month = jul,
  issn = "0024-3795",
  doi = "10.1016/0024-3795(92)90031-5"
  }

@Article{joka02,
  author = "Jonsson, Isak and Kågström, Bo",
  title = "Recursive Blocked Algorithms for Solving Triangular Systems—Part {I}:
           {One}-Sided and Coupled {Sylvester}-Type Matrix Equations",
  journal = j-TOMS,
  year = 2002,
  volume = 28,
  number = 4,
  pages = "392-415",
  month = dec,
  doi = "10.1145/592843.592845"
  }

@Article{lhpw87,
  author = "Laub, A. and Heath, M. and Paige, C. and Ward, R.",
  title = "Computation of System Balancing Transformations and Other Applications
                  of Simultaneous Diagonalization Algorithms",
  journal = j-IEEE-AC,
  year = 1987,
  volume = 32,
  number = 2,
  pages = "115–122",
  month = feb,
  issn = "0018-9286",
  doi = "10.1109/tac.1987.1104549"


}

@Article{lzz20,
  author = "Liu, Zhongyun and Zhou, Yang and Zhang, Yulin",
  title = "On inexact alternating direction implicit iteration for continuous
                  {Sylvester} equations",
  journal = j-NLAA,
  year = 2020,
  volume = 27,
  number = 5,
  month = jul,
  issn = "1099-1506",
  doi = "10.1002/nla.2320"
  }

@Book{saad03,
  author = "Saad, Yousef",
  title = "Iterative Methods for Sparse Linear Systems",
  publisher = pub-SIAM,
  year = 2003,
  address = pub-SIAM:adr,
  edition = "2nd",
  isbn = "0898715342"
  }

@Article{simo16,
  author = "Simoncini, Valeria",
  title = "Computational Methods for Linear Matrix Equations",
  journal = j-SIREV,
  year = 2016,
  volume = 58,
  number = 3,
  pages = "377-441",
  month = jan,
  issn = "1095-7200",
  doi = "10.1137/130912839"
  }

@Article{soan02,
  author = "Sorensen, D. C. and Antoulas, A. C.",
  title = "The {Sylvester} equation and approximate balanced reduction",
  journal = j-LAA,
  year = 2002,
  volume = "351–352",
  pages = "671–700",
  month = aug,
  issn = "0024-3795",
  doi = "10.1016/s0024-3795(02)00283-5"
  }

@Book{svr08,
  author = "Wilhelmus H. A. Schilders and Henk A. Vorst and Joost Rommes",
  title = "Model Order Reduction: Theory, Research Aspects and Applications",
  publisher = pub-Springer,
  year = 2008,
  address = pub-Springer:adr-BH,
  pages = "xi+471",
  doi = "10.1007/978-3-540-78841-6",
  isbn = "978-3-540-78841-6"
  }

@InBook{vand91,
  author = "Van Dooren, Paul M.",
  title = "Structured Linear Algebra Problems in Digital Signal Processing",
  publisher = pub-Springer,
  year = 1991,
  address = pub-Springer:adr-BH,
  pages = "361–384",
  booktitle = "Numerical Linear Algebra, Digital Signal Processing and Parallel
                  Algorithms",
  isbn = 9783642755361,
  doi = "10.1007/978-3-642-75536-1_17"
  }

@Article{vara79,
  author = "Varah, J. M.",
  title = "On the Separation of Two Matrices",
  journal = j-SINUM,
  year = 1979,
  volume = 16,
  number = 2,
  pages = "216-222",
  month = apr,
  issn = "1095-7170",
  doi = "10.1137/0716016",
  }

@Article{wlm13,
  author = "Wang, Xiang and Li, Wen-Wei and Mao, Liang-Zhi",
  title = "On Positive-Definite and Skew-{Hermitian} Splitting Iteration Methods
                  for Continuous {Sylvester} Equation ${AX + XB = C}$",
  journal = j-COMP-MATH-APPL,
  year = 2013,
  volume = 66,
  number = 11,
  pages = "2352-2361",
  month = dec,
  issn = "0898-1221",
  doi = "10.1016/j.camwa.2013.09.011"
  }

@Article{zwt15,
  author = "Zhou, Rong and Wang, Xiang and Tang, Xiao-Bin",
  title = "A Generalization of the {Hermitian} and Skew-{Hermitian} Splitting
                  Iteration Method for Solving {Sylvester} Equations",
  journal = j-APP-MATH-COMP,
  year = 2015,
  volume = 271,
  pages = "609-617",
  month = nov,
  issn = "0096-3003",
  doi = "10.1016/j.amc.2015.09.027"
  }

@article{aabc21,
  author = "Ahmad Abdelfattah and Hartwig Anzt and Erik G. Boman and Erin
            Carson and Terry Cojean and Jack Dongarra and Alyson Fox
            and Mark Gates and Nicholas J. Higham and Xiaoye S. Li and
            Jennifer Loe and Piotr Luszczek and Srikara Pranesh and
            Siva Rajamanickam and Tobias Ribizel and Barry F. Smith
            and Kasia Swirydowicz and Stephen Thomas and Stanimire
            Tomov and Yaohung M. Tsai and Ulrike Meier Yang",
  title = "A Survey of Numerical Linear Algebra Methods Utilizing
           Mixed-Precision Arithmetic",
  journal = j-IJHPCA,
  pages = "344-369",
  volume = 35,
  number = 4,
  month = mar,
  year = 2021,
  doi = "10.1177/10943420211003313",
  created = "2020.07.15",
  updated = "2021.06.20"
        }

@article{cahi17,
  author = "Erin Carson and Nicholas J. Higham",
  title = "A New Analysis of Iterative Refinement and its Application to
           Accurate Solution of Ill-Conditioned Sparse Linear Systems",
  journal = j-SISC,
  volume = 39,
  number = 6,
  pages = "A2834-A2856",
  year = 2017,
  doi = "10.1137/17M1122918",
  created = "2017.03.28",
  updated = "2017.12.11"
       }

@article{cahi18,
  author = "Erin Carson and Nicholas J. Higham",
  title = "Accelerating the Solution of Linear Systems by Iterative
           Refinement in Three Precisions",
  journal = j-SISC,
  volume = 40,
  number = 2,
  pages = "A817-A847",
  year = 2018,
  doi = "10.1137/17M1140819",
  created = "2017.07.22",
  updated = "2018.03.16"
        }

@article{dahi03,
  author = "Philip I. Davies and Nicholas J. Higham",
  title = "A {Schur--Parlett} Algorithm for Computing Matrix Functions",
  journal = j-SIMAX,
  volume = 25,
  number = 2,
  pages = "464-485",
  year = 2003,
  doi = "10.1137/S0895479802410815",
  bibdate = "Thu Sep 18 09:13:12 GMT 2003",
  created = "2003.09.18",
  updated = "2014.05.25"
        }

@article{high93y,
  author = "Nicholas J. Higham",
  title = "Perturbation Theory and Backward Error for {$AX-XB=C$}",
  journal = j-BIT,
  volume = 33,
  pages = "124-136",
  year = 1993,
  doi = "10.1007/978-94-015-8196-7_39",
  created = "1994.01.01",
  updated = "2014.06.20"
        }

@book{high:ASNA2,
  author = "Nicholas J. Higham",
  title = "Accuracy and Stability of Numerical Algorithms",
  publisher = pub-SIAM,
  address =   pub-SIAM:adr,
  edition = "Second",
  year = 2002,
  pages = "xxx+680",
  doi = "10.1137/1.9780898718027",
  isbn = "0-89871-521-0",
  bibdate = "Mon Aug 19 13:30:30 GMT 2002",
  created = "2002.08.19",
  updated = "2014.05.30"
        }

@article{high97i,
  author = "Nicholas J. Higham",
  title = "Iterative Refinement for Linear Systems and {LAPACK}",
  journal = j-IMAJNA,
  volume = 17,
  number = 4,
  pages = "495-509",
  year = 1997,
  doi = "10.1093/imanum/17.4.495",
  bibdate = "Mon Nov 10 11:33:49 GMT 1997",
  created = "1997.11.10",
  updated = "1997.11.10"
        }

@book{high:FM,
  author = "Nicholas J. Higham",
  title = "Functions of Matrices: {Theory} and Computation",
  publisher = pub-SIAM,
  address =   pub-SIAM:adr,
  year = 2008,
  pages = "xx+425",
  doi = "10.1137/1.9780898717778",
  isbn = "978-0-898716-46-7",
  bibdate = "Thu Feb 7 21:51:54 GMT 2008",
  created = "2008.02.07",
  updated = "2014.05.30"
        }

@article{hili21,
  author = "Nicholas J. Higham and Xiaobo Liu",
  title = "A Multiprecision Derivative-Free {Schur--Parlett} Algorithm for
           Computing Matrix Functions",
  journal = j-SIMAX,
  volume = 42,
  number = 3,
  pages = "1401-1422",
  year = 2021,
  doi = "10.1137/20m1365326",
  created = "2020.09.07",
  updated = "2021.08.06"
       }

@article{hima22,
  author = "Nicholas J. Higham and Theo Mary",
  title = "Mixed Precision Algorithms in Numerical Linear Algebra",
  journal = j-AN,
  volume = 31,
  pages = "347-414",
  month = may,
  year = 2022,
  doi = "10.1017/s0962492922000022",
  created = "2021.12.29",
  updated = "2022.06.09"
        }

@phdthesis{kohl21,
    author = "Martin K{\"o}hler",
    title =  "Approximate Solution of Non-Symmetric Generalized Eigenvalue Problems and Linear Matrix Equations on {HPC}-Platforms",
    school = "Logos Verlag, Berlin",
    year = "2021",
    url = {https://pure.mpg.de/pubman/faces/ViewItemOverviewPage.jsp?itemId=item_3326484}
    }

@article{pera55,
  author = "Peaceman, Donald W. and {Rachford, Jr, Henry H.}",
  title = "The Numerical Solution of Parabolic and Elliptic Differential Equations",
  journal = j-SIAM,
  volume = 3,
  number = 1,
  pages = "28-41",
  year = 1955,
  doi = "10.1137/0103003"
        }

@article{luwa91,
  author = "Lu, An and Wachspress, Eugene L.",
  title = "Solution of lyapunov equations by alternating direction implicit iteration",
  journal = j-COMP-MATH-APPL,
  volume = 21,
  number = 9,
  pages = "43-58",
  year = 1991,
  doi = "10.1016/0898-1221(91)90124-M"
        }

@misc{nv25,
  author = "NVIDIA",
  title = "NVIDIA {Blackwell} Architecture Technical Brief",
  type =   "Tech. report",
  year = "2025",
  url = "https://resources.nvidia.com/en-us-blackwell-architecture"
       }
